\documentclass[notitlepage,leqno,10pt]{article}
\usepackage{latexsym}
\usepackage{amsmath}
\usepackage{amsfonts}
\usepackage{amssymb}
\usepackage{amscd}
\renewcommand{\a }{\alpha }
\renewcommand{\b }{\beta }
\renewcommand{\d}{\delta }
\newcommand{\D }{\Delta }
\newcommand{\e }{\varepsilon }
\newcommand{\g }{\gamma}
\renewcommand{\i }{\iota}
\newcommand{\G }{\Gamma }
\renewcommand{\l }{\lambda }
\renewcommand{\L }{\Lambda }
\newcommand{\n }{\nabla }
\newcommand{\s }{\sigma }
\renewcommand{\O }{\Omega }
\renewcommand{\o }{\omega }
\newcommand{\ov}{\overline}
\newcommand{\intbar}{\mathop{\int\makebox(-13.5,0){\rule[4pt]{.7em}{0.3pt}}%
\kern-6pt}\nolimits}
\newcommand{\wtilde }{\widetilde}

\newcommand{\be}{\begin{equation}}
\newcommand{\ee}{\end{equation}}
\newcommand{\bes}{\begin{equation*}}
\newcommand{\ees}{\end{equation*}}
\newcommand{\ba}{\begin{eqnarray}}
\newcommand{\ea}{\end{eqnarray}}
\newcommand{\bas}{\begin{eqnarray*}}
\newcommand{\eas}{\end{eqnarray*}}
\newenvironment{pf}{\noindent{\sc Proof}.\enspace}{\rule{2mm}{2mm}\medskip}

\newtheorem{remark}{Remark}[section]
\newcommand{\R}{\mathbb{R}}

\newcommand{\N}{\mathbb{N}}

\author{ Mohameden  Ahmedou \; $\&$\;\; Mohamed Ben Ayed }

\date{}

\title{\bf  Non simple blow ups for   the Nirenberg problem on half spheres}
\begin{document}

\newtheorem{lem}{Lemma}[section]
\newtheorem{pro}[lem]{Proposition}
\newtheorem{thm}{Theorem}[section]
\newtheorem{cor}[lem]{Corollary}

\maketitle

\centerline{   \emph{dedicated to the  memory of Prof. Louis Nirenberg }}

\begin{center}
{\bf Abstract }
\end{center}
In this paper we study a  Nirenberg type  problem on standard half spheres $(\mathbb{S}^n_+,g_0)$ consisting   of finding conformal metrics of prescribed scalar curvature  and zero boundary mean curvature. This problem amounts to solve the following boundary value problem involving the critical Sobolev exponent:
\begin{equation*}
 (\mathcal{P})  \quad
\begin{cases}
   -\D_{g_0} u \, + \, \frac{n(n-2)}{4} u \,  =  K \, u^{\frac{n+2}{n-2}},\, u > 0 \quad   \mbox{in } \mathbb{S}^n_+, \\
  \frac{\partial u}{\partial \nu }\, =\, 0  \quad  \mbox{on }  \partial \mathbb{S}^n_+,
 \end{cases}
 \end{equation*}
 where $K \in C^3(\mathbb{S}^n_+)$  is  a positive  function.
 We construct, under generic conditions on the function $K$,  finite energy solutions of a subcritical approximation of $(\mathcal{P})$ on half spheres of dimension $n \geq 5$, which exhibit multiple blow up of \emph{cluster-type} at the same  boundary point. These solutions may have zero or non zero weak limit and may  develop clusters at different boundary points.  Such a blow up phenomena   on half spheres drastically contrast  with the  case of the Nirenberg problem on spheres, where non simple blow up for finite energy solutions cannot occur and unveils  an unexpected connection with vortex type problems arising in Euler equations in fluid dynamic and mean fields type equations in mathematical physics.
 We construct also, under suitable conditions on the restriction of $K$ on $\partial \mathbb{S}^n_+$,  approximate solutions of arbitrarily large energy and Morse index.

\begin{center}


\centerline{\bf AMS subject classification:  35A01, 58J05, 58E05.}

\end{center}

\tableofcontents

\section{Introduction and statement of the results}

For the $n-$dimensional  sphere  $ \mathbb{S}^n, \, n \geq 3$ endowed with its standard metric $g_0$ and a given  positive function $K \in C^2(\mathbb{S}^n)$,   Nirenberg asked in the early seventieth  the following question:
 Can the function  $K$  be realized as the scalar curvature of a metric $g$, conformally equivalent to $g_0$? \\
  Writing the conformal metric as  $g:= u^{{4}/{(n-2)}} g_0$, the above question amounts to solve the following nonlinear elliptic equation involving the Sobolev critical exponent:
\begin{equation}\label{eq:Sn}
 (\mathcal{NP}) \quad  -\D_{g_0} u \, + \, \frac{n(n-2)}{4} u \,  =  K \, u^{\frac{n+2}{n-2}}, \quad \, u > 0\,  \mbox{ in }\,  \mathbb{S}^n.
\end{equation}
\noindent
This problem has a variational structure, however  the associated Euler-Lagrange  functional does not satisfy the Palais-Smale condition due to the critical growth of the nonlinearity. \\
To overcome such a difficulty, one can  lower   the  critical exponent by considering the following subcritical  approximation of the Nirenberg problem $(\mathcal{NP})$:

   \begin{equation}\label{eq:ep}
  (\mathcal{NP}_{\e}) \quad  -\D_{g_0} u \, + \, \frac{n(n-2)}{4} u \,  =  K \, u^{\frac{n+2}{n-2} - \e}, \quad  \, u > 0 \, \mbox{ in }\,  \mathbb{S}^n,
\end{equation}
where $\e > 0$ is a small parameter.
In this way one recovers the compactness and  then studies the behavior of a solution $u_{\e}$ of $(\mathcal{NP}_{\e})$ as the parameter $\e$ goes to zero. Thanks to  elliptic estimates there are two alternatives: either $||u_{\e}||_{L^{\infty}}$ remains uniformly bounded and in this case the solutions $u_{\e}$ converge in the $C^{2,\a}-$topology  as $\e \to 0$ to a solution of $(\mathcal{NP})$ or the sequence of solutions $u_{\e}$ blows up. In this latter case, assuming that $u_\e$ is a family of energy-bounded solutions,  it follows from the concentration compactness principle that $ u_{\e}^{{2n}/{(n-2)}}  \mathcal{L}^n$, where $\mathcal{L}^n$ denotes the $n$-dimensional Lesbegue measure,  converges in the sense of measures to a sum of Dirac masses, see \cite{Lions, Struwe}. Actually a  refined blow up analysis of such blowing up solutions  has been initiated   by R. Schoen \cite{Schoen, SZ, KMS09}  and developed by  Y.Y. Li \cite{yyli1, yyli2},  C.C. Chen and  C.S Lin \cite{CL1, CL2}) and Druet-Hebey-Robert \cite{DHR}. It follows from such a blow up analysis that solutions concentrate at critical points of $K$. Moreover, in dimension $n=3,4$ and  under the assumption that  $\D K$ does not vanish at any critical point, it has been proved that all blow up points are \emph{isolated simple}, that is locally the blowing up solution has the energy of \emph{one bubble}, see please definition 0.3 in \cite{yyli1}. Moreover concentration  occurs only at critical points of $K$ with $\D K < 0$. This has been proved by Yanyan Li   \cite{yyli1, yyli2} for the dimensions $n=3,4$. Furthermore on  spheres of  dimensions $n \geq 5$  and  under the additional assumption that the energy of the solutions is uniformly bounded,    A. Malchiodi and M. Mayer   \cite{MM} proved that all blow up points are isolated simple. See also A. Bahri \cite{Bahri-Invariant} where a corresponding property has been proved for \emph{critical points at infinity} of the associated variational problem. Furthermore the assumption that the energy  is bounded is necessary to rule out \emph{non simple blow up} as shown by C.C. Chen and C.S. Lin \cite{CL}.

In this paper we consider a Nirenberg type  problem on standard half spheres $(\mathbb{S}^n_+,g_0)$ consisting of  prescribing  simultaneously   the scalar curvature to be a function $K \in C^3(\mathbb{S}^n_+)$ and the boundary mean curvature to be zero. This amounts to solve the following boundary value problem
\begin{equation}
(\mathcal{P}) \quad
\begin{cases}
   -\D_{g_0} u \, + \, \frac{n(n-2)}{4} u \,  =  K \, u^{\frac{n+2}{n-2}},\,  u > 0 & \mbox{ in } \mathbb{S}^n_+, \\
  \frac{\partial u}{\partial \nu }\, =\, 0  & \mbox{on } \partial \mathbb{S}^n_+.
\end{cases}
\end{equation}
This problem has been studied  on half spheres of dimensions $n=3,4$. See  the papers \cite{yyli, DMOA, BEOA, BEO, BGO} and the references therein.\\
Here also in order to recover compactness  one considers the  following subcritical approximation

 \begin{equation}
(\mathcal{P}_{\e}) \quad
\begin{cases}
   -\D_{g_0} u \, + \, \frac{n(n-2)}{4} u \,  =  K \, u^{\frac{n+2}{n-2} - \e}, \,  u > 0 & \mbox{ in } \mathbb{S}^n_+, \\
  \frac{\partial u}{\partial \nu }\, =\, 0  & \mbox{on } \partial \mathbb{S}^n_+.
\end{cases}
\end{equation}
Just as above, there are two alternatives for the behavior of a sequence of solutions $u_{\e}$ of $(\mathcal{P}_{\e})$. Either the $||u_{\e}||_{L^{\infty}}$ remains uniformly bounded or it blows up and if it does, assuming that $u_\e$ is a family of energy-bounded solutions, then it follows that $u^{{2n}/{(n-2)}}_{\e}  \mathcal{L}^n$ converges to a sum of Dirac masses, some of them are in the interior and the others are located on the boundary. The interior points are  critical points of $K$  satisfying  that $\D K \leq 0$ and the boundary points are critical points of $K_1$ the restriction of $K$ on the boundary  satisfying that $\partial_{\nu} K \geq  0$. See \cite{DMOA}  and \cite{BEO}.
Moreover a refined blow up analysis, under the assumption that $\D K \neq 0$ at interior critical points of $K$ and that $\partial_{\nu} K \neq 0$ at critical points of $K_1$,  shows that in the dimension $n=3$  all blow up points are isolated simple, see \cite{yyli} and  \cite{DMOA}. Furthermore,  under additional condition on $K_1$,  it has been proved in \cite{BGO} that on four dimensional half spheres  all blow up points are isolated simple.\\
On half spheres of dimension $n \geq 5$ the problem $(\mathcal{P}_{\e})$ admits  blowing up solutions having only    \emph{ isolated simple blow points}. Indeed we have:
\begin{thm} \label{t:13}
Let $K$ be a positive smooth function on $\ov{\mathbb{S}^n_+}$  with $n \geq 5$ and let $z_1,\cdots,z_m \in \partial \mathbb{S}^n_+$ be  non degenerate critical points of $K_1:= K_{\lfloor \partial \mathbb{S}^n_+}$ which satisfy that $(\partial_{\nu} K)(z_i) > 0$ for each $i$ and let $y_{m+1},\cdots, y_{m+\ell}$ be non degenerate critical points of $K$ with $\D K(y_i) < 0$ for each $i > m$.\\
Then there exists a  sequence of solutions $u_{\e}$ of $(\mathcal{P}_\e)$  which converges weakly to $0$ and blows up at $z_1, \cdots, z_m, y_{m+1}, \cdots,y_{m+\ell}$ and for each $ i $, it holds
$$
 \lim_{r \to 0}  \lim_{ \e  \to 0} \int_{B_r(z_i)} K u_{\e}^{\frac{2n}{n-2}}  \, =  \, \underbrace{\frac{S_n }{K(z_i)^{(n-2)/2}}}_\text{Energy of  1 boundary bubble}  \mbox{ and }  \lim_{r \to 0}  \lim_{ \e  \to 0} \int_{B_r(y_i)} K u_{\e}^{\frac{2n}{n-2}} \, = \, \underbrace{ \frac{2\, S_n }{K(y_i)^{(n-2)/2}}}_\text{Energy of 1 interior bubble}
$$
where $B_r(z)$ is a geodesic ball and the universal constant $S_n:= c_0^{2n/(n-2)}\int_{\R^n_+} \frac{1}{(1+ |y|^2)^n}$, where $c_0:= (n(n-2))^{(n-2)/4}$,  is the energy of one \emph{boundary bubble} (for $K\equiv 1$).
Hence    each blow up point is an isolated simple blow up one.
\end{thm}

\begin{remark}
\begin{enumerate}
  \item The above theorem holds  true on three and four dimensional half spheres. See \cite{DMOA, BGO}.
  \item  The above theorem gives the existence of  blowing up solutions having only interior blow up points (by taking $m = 0$ ) or  only boundary blow up points (by taking $\ell = 0$).
\end{enumerate}
\end{remark}

Before resuming our investigation of the blow up phenomena we recall the definition of non degenerate solution of $(\mathcal{P})$:
\emph{A solution $\o$ of $(\mathcal{P})$ is said to be non degenerate if the linearized operator
\begin{equation}\label{eq:lin}
\mathcal{L}_{\omega}(\varphi):= - \D_{g_{0} } \varphi \, + \frac{n(n-2)}{4} \varphi \, - \, \frac{n+2}{n-2} K \o^{\frac{4}{n-2}} \varphi
\end{equation}
does not admit  zero as an eigenvalue.}\\
Next we want to consider the question of existence of  blowing up solutions which are close to a combination of  a sum of bubbles and a solution of $(\mathcal{P})$.
Namely  we prove that  there are   blowing up solutions of $(\mathcal{P}_{\e})$ with non zero weak limit having interior as well as boundary blow ups which are all \emph{isolated simple}. More precisely our result can be stated as follows:

\begin{thm} \label{t:18}
Let $K$ be a positive smooth function on $\ov{\mathbb{S}^n_+}$,  $z_1,\cdots,z_m \in \partial \mathbb{S}^n_+$ be  non degenerate critical points of $K_1:= K_{\lfloor \partial \mathbb{S}^n_+}$ which satisfy that $(\partial_{\nu} K)(z_i) > 0$ for each $i$ ($m\in \N\cup \{0\}$) and  $y_{m+1},\cdots, y_{m+\ell}$ be non degenerate critical points of $K$ with $\D K(y_i) < 0$ for each $i > m$ ( $\ell \in \N\cup \{0\}$). Furthermore let $\o$ be a non degenerate solution of $ (\mathcal{P})$.
   \begin{itemize}
     \item[(i)]  If $n \geq 5 $
     there exists a  sequence of solutions $u_{\e}$ of $(\mathcal{P}_\e)$  which converges weakly to $\o$ and blows up at $(z_1, \cdots, z_m)$  (with $m\geq 1$) and for each $ i $, it holds
$$
 \lim_{r \to 0}   \lim_{ \e  \to 0} \int_{B_r(z_i)} K u_{\e}^{({2n}/{(n-2)}) } \, =  \, \frac{S_n }{K(z_i)^{(n-2)/2}}.
$$
     \item[(ii)]  If $n \geq 7$
     there exists a  sequence of solutions $u_{\e}$ of $(\mathcal{P}_\e)$  which converges weakly to $\o$ and blows up at $z_1, \cdots, z_m, y_{m+1}, \cdots,y_{m+\ell}$  (with $m \geq 0$ and $\ell \geq 1$) and for each $ i $, it holds
$$
 \lim_{r \to 0}  \lim_{ \e  \to 0} \int_{B_r(z_i)} K u_{\e}^{\frac{2n}{n-2} } \, =  \, \frac{S_n }{K(z_i)^{(n-2)/2}} \quad  \mbox{ and } \quad \lim_{r \to 0}     \lim_{ \e  \to 0} \int_{B_r(y_i)} K u_{\e}^{\frac{2n}{n-2} } \, = \, \frac{ 2 \,S_n }{K(y_i)^{(n-2)/2}} .
$$
  \end{itemize}
All the blow up points are isolated simple. \\
   We call such a blow up behavior, \emph{blow up with residual mass}.
\end{thm}

\begin{remark}
The existence of blow up with residual mass in dimension $n =6$ is an open question even for the case of closed spheres. See please the remarks  ending the proof of Theorem \ref{t:15} for the analytical  features behind the failure of our argument in this dimension.
\end{remark}

In contrast to  the above  results and also to the case of  closed spheres, we construct in this paper a sequence of finite energy solutions of $(\mathcal{P}_{\e})$ on half spheres of dimension $n \geq 5$, which exhibit multiple blow up at the same  boundary point.

Indeed it turns out that the existence of solutions blowing up at the same boundary point with the energy of $m$ boundary bubbles, $ m \geq 2$ is related to the existence of critical points of   the following Kirchhoff-Routh type  function:
\be \label{dF2} \mathcal{F}_{z,m}:  \mathbb{F}_m( T_z (\partial \mathbb{S}^n_+))\to \R ; \, \mathcal{F}_{z,m}(\xi_1,\cdots,\xi_m):= \frac{1}{2}\sum_{i=1}^{m} D^2K_1(z)( \xi_i , \xi_i )  + \sum_{1 \leq i < j \leq m}\frac{1}{| \xi_i - \xi_j | ^{n-2}}  .\ee
where $ \mathbb{F}_m( T_z (\partial \mathbb{S}^n_+)):=  \{   (\xi_1,\cdots,\xi_m); \, \xi_i \neq \xi_j  \in T_z (\partial \mathbb{S}^n_+) \, \, \mbox{ for } i \neq j  \} $.\\
We point out that similar   hamiltonian type functions appear in the characterization of the location of the concentration points in mean field type equations in mathematical physics  and Euler equation in fluid dynamic. The relevance of such vortex  type problems  in Yamabe  type equations has been discovered by  Thizy-V\'etois \cite{Thizy-Vetois} and Pistoia-Vaira \cite{Pistoia-Vaira}. \\
Next we prove that every non-degenerate critical point of the functional $ \mathcal{F}_{z,m}$, where $z$ is a boundary point satisfying the assumption of Theorem \ref{t:01} gives rise to bubbling off solutions of high energy and morse index. We recall that the \emph{morse index} $Morse(\o)$ of a solution $\o$ of $(\mathcal{P})$ is the dimension of the space of negativity of the bilinear form associated to the linearized operator $\mathcal{L}_{\o}$ (defined in \eqref{eq:lin}).
Namely we have

\begin{thm} \label{t:11}
Let $K$ be a positive smooth function on $\ov{\mathbb{S}^n_+}$  with $n \geq 5$ and $z \in \partial \mathbb{S}^n_+$ be a non degenerate critical point of $K_1:= K_{\lfloor \partial \mathbb{S}^n_+}$ satisfying  $(\partial_{\nu} K)(z) > 0$.\\
For $m \geq 2$, assume that the function $  \mathcal{F}_{z,m}$ has a non degenerate critical point. Then there exists a  sequence of solutions $u_{\e}$ of $(\mathcal{P}_\e)$ having the energy of $m$ boundary bubbles, which blows up at $z$. That is
$$
  \lim_{r \to 0}  \lim_{ \e  \to 0} \int_{B_r(z)} K u_{\e}^{({2n}/{(n-2)}) } \, = \,\underbrace{ m \, {\frac{1}{K(z)^{(n-2)/2}} S_n}}_{\text{Energy of m  boundary bubbles}} .
$$
Furthermore, the morse index of the solution $u_\e$ is lower bounded by $m$.\\
We call such a blow up point  non simple of order m.
\end{thm}

\begin{remark}\label{rem2}
\begin{itemize}
  \item
  The above non simple blow up is of cluster type and we point out that non simple blow up of cluster type have been recently constructed for  some perturbations  of   some Yamabe  type equations  on manifolds which are not locally  conformally flat  by  Thizy-V\'etois \cite{Thizy-Vetois} and Pistoia-Vaira \cite{Pistoia-Vaira}.
  \item
We point out that assuming the function   $\mathcal{F}_{z,m}$ has a critical point implies that  the  boundary point $z$,  where the concentration occurs, is not  a local maxima. See  statement $(ii) $ in Proposition \ref{pF}.
\item
A non simple blow up phenomenon has been proved by P.Esposito and J.Wei \cite{Esposito-Wei} for the Sinh-Gordon equation under Neumann boundary conditions.
\end{itemize}
\end{remark}

As an immediate corollary of Theorem \ref{t:11} we have the following result

\begin{cor} \label{t:01} Let $K$ be a positive smooth function on $\ov{\mathbb{S}^n_+}$  with $n \geq 5$ and let $z$ be a non-degenerate critical point  of $K_1:= K_{\lfloor \partial \mathbb{S}^n_+}$ which  satisfies that  $\partial_{\nu} K(z) > 0$.\\  Suppose that  $D^2K_1(z)$ has at least one  simple positive eigenvalue. Then  there exists a sequence of  solutions $u_{\e}$ of $(\mathcal{P}_\e)$ having the energy of  2 boundary bubbles and which blows up at the boundary point  $z$. That is $$  \lim_{r \to 0}  \lim_{ \e  \to 0} \int_{B_r(z)} K u_{\e}^{({2n}/{(n-2)})} \, = \, 2 \, \underbrace{\frac{1}{K(z)^{(n-2)/2}} S_n}_\text{Energy of 1 boundary bubble}. $$ \end{cor}

 Furthermore, as a corollary of Theorem \ref{t:11} and Proposition \ref{p:dF234} we obtain the  existence result of  blowing up solutions of $ (\mathcal{P}_{\e})$ with arbitrarily large energy and Morse index.  More precisely we prove the following theorem:

\begin{cor}\label{t:inftyenergy}
Let $n \geq 5$ and assume that the function $K_1:= K_{\lfloor{\partial \mathbb{S}^n_+}}$ has a non degenerate critical point $z \in \partial \mathbb{S}^n_+$ such that $\partial_{\nu}K(z)  >0$  and the hessian  matrix  $D^2 K_1(z)$ has only one simple positive eigenvalue and the others are negative. Then for every $m \in \N$ there exists a solution $u_{\e,m}$ of $(\mathcal{P}_{\e})$ such that
$$
  \lim_{m \to \infty}  \lim_{ \e \to 0} \int_{\mathbb{S}^n_+} K u_{\e,m}^{\frac{2n}{n-2}} \, =  \, + \infty \quad \mbox{ and }  \quad  Morse(u_{\e,m}) \to \infty, \mbox{ as } m \to \infty .$$
\end{cor}

In the next result we generalize Theorem \ref{t:11} by constructing a solution of $(\mathcal{P}_{\e})$ having  \emph{clusters} of non simple blow up points  at different boundary points.
 Furthermore, the solution can have  also interior blow up points. Namely we prove:

\begin{thm} \label{t:12}
Let $K$ be a positive smooth function on $\ov{\mathbb{S}^n_+}$  with $n \geq 5$ and let $z_1,\cdots,z_m \in \partial \mathbb{S}^n_+$  (with $m \geq 0$) be  non degenerate critical points of $K_1:= K_{\lfloor \partial \mathbb{S}^n_+}$, which  satisfy that $(\partial_{\nu} K)(z_i) > 0$ for each $i$  and  $y_{m+1},\cdots, y_{m+\ell}$ be non degenerate critical points of $K$ with $\D K(y_i) < 0$ for each $i > m$ ( $\ell \geq 0$).\\
For $1 \leq i \leq m$ let $q_i \geq 1$ and suppose that if $q_i \geq 2$ then  the function $  \mathcal{F}_{z_i,q_i}$ has a non degenerate critical point. Then there exists a  sequence of solutions $u_{\e}$ of $(\mathcal{P}_\e)$  which  converges weakly to $0$ and  blows up at $z_1, \cdots, z_m$,  $y_{m+1},\cdots, y_{m+\ell}$. Precisely, for each $i$, it holds
$$
  \lim_{r \to 0}   \lim_{ \e  \to 0} \int_{B_r(z_i)} K u_{\e}^{({2n}/{(n-2)}) } \, = \, q_i \, \frac{S_n}{K(z_i)^{(n-2)/2}}   \mbox{ and }  \lim_{r \to 0}   \lim_{ \e  \to 0} \int_{B_r(y_i)} K u_{\e}^{({2n}/{(n-2)}) } \, = \,  \frac{ 2\,  S_n }{K(y_i)^{(n-2)/2}} .
$$
\end{thm}

  We notice that, in the case where $\ell = 0$ and $m \geq 1$, there are  only boundary blow up points and, in this case, the second limit vanishes, while for  $m=0$ and $\ell \geq 1$  there are  only interior blow up points and in this case the first limit vanishes.\\
We also point out that  Theorem \ref{t:12}  combined with the bubbling off analysis in \cite{AB20b}, describes the full blow up picture  for bounded energy  solutions converging  weakly to $0$. Indeed, if $q_i = 1$ for each $i$ the statement reduces to the one of   Theorem \ref{t:13}, while  if $m=1$, $q_1\geq 2$ and $\ell =0$, one recovers  Theorems \ref{t:01} and \ref{t:11}. We recall that for the closed case, similar  results have been obtained by Malchiodi-Mayer in \cite{MM-2,MM}.

\medspace

In the next theorem we construct  blowing up solutions of $(\mathcal{P}_{\e})$ having non zero weak limit and exhibiting a blow up of cluster type at a  boundary point. Namely we prove:

\begin{thm} \label{t:15}
Let $\o$ be a non degenerate solution of $(\mathcal{P})$.
\begin{enumerate}
\item[a)] For $n \geq 5$,  under the same assumptions of Theorem \ref{t:12} with $\ell = 0$ and $m \geq 1$ (that is we have only boundary blow up points), there exists a  sequence of solutions $u_{\e}$ of $(\mathcal{P}_\e)$  which converges weakly to $\o$  with the same properties than Theorem \ref{t:12}.
\item[b)] Let $ n \geq 7$. Under the same assumptions of Theorem \ref{t:12} with $\ell \geq 1$ and $m \geq 0$ (that is we have at least one interior blow up point), there exists a  sequence of solutions $u_{\e}$ of $(\mathcal{P}_\e)$  which converges weakly to $\o$  with the same properties than Theorem \ref{t:12}.
\end{enumerate}
\end{thm}

Regarding the method of the proof of the above theorems some remarks are in order. Although we use  a finite dimensional framework similar to the one used  by Bahri-Li-Rey \cite{BLR} to deal with Yamabe type equations on domains, there is a major difference  in the analytical feature of  the two  frameworks. Indeed Bahri-Li-Rey used in crucial way the refined blow up analysis of R. Schoen \cite{Schoen, KMS09} which proves that all blow up points  of the Yamabe equation are \emph{ isolated simple}. Such a fact simplifies  in a drastic way the asymptotic expansion of the gradient near potential blow up points and hence the construction of blowing up solutions. In particular it implies that the  mutual interaction of two bubbles  $\d_{a_i,\l_i}$  and $\d_{a_j,\l_j}$ is of the order of ${G(a_i,a_j)}/{(\l_i\l_j)^{{(n-2)}/{2}}}$, where $G(.,.)$ denotes the Green's function. In contrast with the situation considered by Bahri-Li-Rey \cite{BLR} we are in a situation where the points are very close to each other and it is very challenging to compute the leading term in the mutual interaction. In fact with respect to the cases of the Nirenberg problem  on spheres and  of the Yamabe problem, the main analytical feature  relies in that fact that at every critical point $z \in \partial \mathbb{S}^n_+$ of  the restriction of $K$ on the boundary $K_1:= K_{\lfloor{\partial \mathbb{S}^n_+}}$, which is not a local maximum and satisfying that $\frac{\partial K}{\partial \nu}(z) > 0$,    a \emph{ balancing phenomenon } between the self interaction and the mutual interaction of bubbles sitting  near of $z$ occurs. This balancing phenomenon   is the key point of our construction of multibubbling solutions of the approximated problem $(\mathcal{P}_{\e})$ accumulating at the same boundary point $z$. Furthermore the  goals of the two papers are different. Indeed the main result of Bahri-Li-Rey is  the computation of the difference of topology induced by \emph{simple blow points} while our goal in this paper is to construct \emph{non simple blow up points.} In  a subsequent paper, we will compute  the difference of topology induced by  these \emph{non simple blow up points} and show that the computation is much more involved compared with the  \emph{simple blow up} case. 

\medspace

The sequel  of this paper is organized as follows: in Section 2 we  set up the variational framework associated to equation $(\mathcal{P}_{\e})$, recall its related Euler-Lagrange functional  and its finite dimensional reduction. In Section 3 we expand the gradient of the Euler-Lagrange functional in the neighborhood of highly concentrated \emph{bubbles}. Section 4 is devoted to the study of the critical point of a Kirchhoff-Routh type functional while Section 5 is devoted to prove our main results. Finally we collect in the appendix useful estimates of the \emph{standard bubble}.

\section{ Variational  framework and preliminaries}

In this section we set up the general framework and introduce some notation.\\
We start by recalling  the variational framework  of  the  boundary value problem $(\mathcal{P}_{\e})$. Namely we recall that  solutions of $(\mathcal{P}_{\e})$ are in one to one correspondence with the positive critical points of the functional
$$ I_\e(u)=\frac{1}{2}\int_{\mathbb{S}^n_+}|\n u|^2+ \frac{n(n-2)}{8} \int_{\mathbb{S}^n_+} u^2-\frac{1}{p+1-\e}\int_{\mathbb{S}^n_+}K|u|^{p+1-\e},\, \,  u \in H^1(\mathbb{S}^n_+) \mbox{ and } p:=\frac{n+2}{n-2}. $$
For convenience matter,  we  perform a stereographic projection  to reduce our problem to $\R^n_+$. Let $\mathcal{D}^{1,2}(\R^n_+)$ denote the completion of $C^\infty_c (\overline{\R^n_+})$ with respect to Dirichlet norm. The stereographic projection $\pi _a$ through a point $a \in \partial \mathbb{S}^n_+$ induces an isometry $\i : H^1(\mathbb{S}^n_+) \to \mathcal{D}^{1,2}(\R^n_+) $ according to the following formula
\begin{align}\label{i44} (\i v)(x)= \Big(\frac{2}{1+|x|^2}\Big)^{(n-2)/{2}}v(\pi _a^{-1}(x)), \qquad v\in H^1(\mathbb{S}^n_+), \, x\in \R^n_+. \end{align}
In particular we have that  for every $v\in H^1(\mathbb{S}^n_+)$
\be \label{i444} \int _{\mathbb{S}^n_+}(|\n v|^2 + \frac{n(n-2)}{4} v^2) = \int _{\R^n_+}|\n (\i v)|^2 \qquad \mbox{and } \qquad\int _{\mathbb{S}^n_+}|v|^{p+1}= \int _{\R^n_+}|\i v|^{p+1}. \ee
In the sequel, we will identify the function $K$ and its composition with the stereographic projection $\pi _a$. We will also identify a point $b$ of $\mathbb{S}^n_+$ and its image by $\pi _a$.
We define
\begin{equation} \label{dal}
\d_{a,\l}(x) := c_0 \frac{\l^{n-2/2}}{ (  \l^2 + 1 + (1-\l^2) \cos d(a,x))^{n-2/2}}  \qquad \mbox{ with } c_0:= (n(n-2))^{(n-2)/4}
\end{equation}
where $d$ is the geodesic distance on  $\ov{\mathbb{S}^n_+}$ and the  constant $c_0$ is  chosen such that $$
- \D \d_{a,\l} \, + \, \frac{n(n-2)}{4} \d_{a,\l} \, = \, \d_{a,\l}^{{(n+2)}/{(n-2)}} \quad \mbox{ in } \mathbb{S}^n_+.
$$
We notice that this function satisfies  $ c / \l^{(n-2)/2} \leq  \d_{a,\l} \leq c  \l^{(n-2)/2} $   uniformly in $\mathbb{S}^n_+$, which implies that
\be\label{lower} | \d_{a,\l} ^{-\e} - 1 | \leq c\,  \e \ln \l \qquad \mbox{ if } \e\ln \l \mbox{ is small}.\ee
Note that, in some computations, we need a precise estimate of $\d_{a,\l} ^{-\e}$ which is given in Lemma \ref{lowerL2}.\\
Moreover, it is easy to see that, using \eqref{i44} with $\pi_{-a}$, the function $\i\delta_{(a,\l)}$ is equal to
\begin{equation}\label{eq:delta}
  \i\delta_{a,\l}(x) = c_0\frac{\l ^{(n-2)/2}}{(1+\l^2|x|^2)^{(n-2)/2}} \, ; \quad x \in \R^n_+
\end{equation}
(see the proof of Lemma \ref{lowerL2} for the details of this change of variables).

Note that, if we use $\pi_b$ (that is with $b \in \partial  \mathbb{S}^n_+$ instead of $-a$),
the function $\i \d_{a,\l}$ will be $$ c_0\frac{\mu ^{(n-2)/2}}{(1+\mu^2| x - \tilde{a} |^2)^{(n-2)/2}} \qquad  \mbox{ with } \quad \begin{cases} \tilde{a} := \frac{(\l^2 -1) \mbox{Proj}_{\R^n} a }{2+(\l^2-1)(1-\cos d(a,b))} \\
 \mu := \frac{2+(\l^2-1)(1-\cos d(a,b))}{2 \l} \end{cases}
$$
where $\mbox{Proj}_{\R^n}$ denotes the projection onto $\R^n$ (for the proof, see page 14 of \cite{BaBr}). In the following we will use $\d_{a,\l}$ even if we are in $\R^n$. These notation  will be assumed to be  understood in the sequel.

Notice that, if the point $a\in \mathbb{S}^n_+$ then the normal derivative on the boundary of the function $\d_{a,\l}$ is not $0$. Hence we need to modify this function as follows:\\
For $a \in \ov{ \mathbb{S}^n_+},$ we define \emph{projected bubble} $\varphi_{a,\l}$ to be the unique solution of
$$
- \D \varphi_{a,\l} \, + \, \frac{n(n-2)}{4} \varphi_{a,\l} \, = \, \d_{a,\l}^{{(n+2)}/{(n-2)}} \quad \mbox{ in }\,  \mathbb{S}^n_+; \quad \frac{\partial \varphi_{a,\l}}{\partial \nu} \, = 0 \mbox{ on } \partial \mathbb{S}^n_+.
$$
We point out that $\varphi_{a,\l} = \d_{a,\l}$ if $a \in  \partial \mathbb{S}^n_+$.

\medskip
 For $\e = 0$, the functional $I_0$ does not satisfy the Palais-Smale condition $(PS)$ and non converging $(PS)$-sequences of positive functions belong to the   set $ V(q,\ell,\tau)$ if their weak limit is zero and in $ V(\o, q,\ell,\tau)$ if their weak limit is $\o$. These sets called in the sequel \emph{neighborhood at infinity} are defined as follows:  \\
 For $q, \ell \in \N \cup \{0\} $ (with $q+\ell \geq 1$) and $\tau > 0$ small
 \begin{align}
  V(q,\ell,\tau):= & \{ u \in H^1(\mathbb{S}^n_+): \, \exists \, \l_i > \tau^{-1}, \,  \exists \,   a_i\in \ov{\mathbb{S}^n_+} \, \mbox{ for } i = 1 \cdots ,q+\ell \, \, \mbox{ s.t. } \, \l_i d(a_i, \partial \mathbb{S}^n_+)  < \tau \\
  &  \mbox{ for } i \leq q \, ;  \, \,   \l_i d(a_i, \partial \mathbb{S}^n_+) > \tau^{-1}   \mbox{ for }  i > q \,  \mbox{ and } \| u - \sum K(a_i)^{(2-n)/4} \varphi_{a_i,\l_i} \| < \tau\}, \nonumber \\
 V(\o, q,\ell,\tau):= & \{ u \in H^1(\mathbb{S}^n_+): u - \o \in  V(q,\ell,\tau) \},
\end{align}
 where $\o$ is a positive solution of $(\mathcal{P})$ and where $\e_{ij}$ denotes  the mutual  interaction between two different bubbles and it is defined by
 \begin{equation}\label{epsilon}
 \e_{ij}:= [\l_i/\l_j + \l_j/ \l_i + ({1}/{2}) \l_i\l_j (1-\cos d( a_i , a_j ) )]^{(2-n)/2} .
 \end{equation}

 \subsection{ Finite dimensional reduction in  the zero weak limit case}

In the case  where the weak limit of non converging  $(PS)$-sequences of positive functions is zero, we proceed as follows to parameterize the neighborhood at infinity $ V( q,\ell,\tau)$. Namely we have the following proposition whose proof  is, up to minor modifications,  identical to the proof of Proposition 7 in  \cite{BCd}.
 \begin{pro} Let $u \in V(q,\ell, \tau)$, then the following minimization problem
 $$ \min _{ \a_i > 0;\l_i > 0; a_i\in \partial  \mathbb{S}^n_+ (i \leq q); a_i \in \mathbb{S}^n_+ ( i > q)} \| u - \sum_{ 1 \leq i \leq  q + \ell } \a_i \varphi_{a_i, \l_i} \| $$
 has a unique solution (up to permutations on the indices).
 Hence, each function $u \in V(q,\ell,\tau)$ can be written as
 \be u =   \sum_{ 1 \leq i \leq q+\ell} \a_i \varphi_{a_i, \l_i}  + v \ee
 where the parameters $\a_i, a_i, \l_i$ are the solution of the previous minimization problem satisfying
 \be\label{alp5}
| \a_i^{4/(n-2)} K(a_i) - 1 | < c \, \tau \, \, \,  \forall \, \, 1 \leq i \leq q+\ell \quad  ; \quad   \l_i d(a_i, \partial \mathbb{S}^n_+) > c \, \tau^{-1}  \, \, \, \forall i > q\ee
  and $v$ satisfies
 \begin{align}\label{V0} & \| v \| \leq \tau \mbox{ and } \langle v, \psi \rangle = 0 \, \,  \mbox{ for each }\\
 & \psi \in \{ \varphi_i ;  \partial \varphi_i / \partial \l_i ;  \partial \varphi_i / \partial a_i^k ,\, \,   i \geq q+1 ,\,  k \leq n\} \cup \{ \d_i ;  \partial \d_i / \partial \l_i ;  \partial \d_i / \partial a_i^k ,\, \,  i \leq q ,\,  k \leq n-1\} .\nonumber \end{align}
 \end{pro}

Regarding the infinite dimensional part $v$ we have the following estimate:
  \begin{pro} \label{pv}
  Let $u:= \sum_{i=1}^{q+\ell} \a_i \varphi_{a_i,\l_i}  \in V(q,\ell, \tau)$ with $\e \ln \l_i$ is small for each $i$. Then there exists a unique function $\ov{v}:= \ov{v}_{\e, \a_i, a_i,\l_i}$ such that
  \be\label{AAAZZZ} \langle \n I_\e (u + \ov{v}) , v \rangle = 0 \quad  \mbox{ for each } v \mbox{ satisfying  \eqref{V0}}.\ee
   Furthermore, it holds that :  $\|\overline{v} \| \leq c R(\e,a,\l)$ where
$$R(\e,a,\l) :=   \e +   \sum_{i=1}^{q+p}  \Big( \frac{| \n K(a_i) | }{\l _i} + \frac{1}{\l _i^2}  \Big) +  \, \begin{cases}
 \sum \e _{ij }^{\frac {n+2}{2(n-2)}}(\ln \e _{ij}^{-1})^{\frac{n+2}{2n}} + \sum_{i > q }\frac{\ln(\l_i d_i)}{(\l_i d_i)^{(n+2)/2}}\mbox{ if } n \geq 6,\\
 \sum \e _{ij }(\ln \e _{ij}^{-1})^{{3}/{5}} + \sum_{i > q} \frac{1}{(\l_i d_i)^{3}}\mbox{ if } n = 5. \end{cases}
$$
   \end{pro}
 The  above proposition   follows using similar arguments as in  the proof of Proposition \ref{wpv} below.

  \subsection{ The non zero weak limit case}

In this subsection we deal with  non converging $PS$-sequences having non zero weak limit. Let $\o$ be a non degenerate solution of $(\mathcal{P})$. Arguing as above we have the following parametrization:
  \begin{pro} Let $u \in V(\o, q,\ell, \tau)$, then the following minimization problem
 $$ \min _{ \a_i > 0;\l_i > 0; a_i\in \partial  \mathbb{S}^n_+ (i \leq q); a_i \in \mathbb{S}^n_+ ( i > q) } \| u - \a_0\o - \sum_{ 1 \leq i \leq q+\ell} \a_i \varphi_{a_i, \l_i} \| $$
 has a unique solution (up to permutations on the indices).\\
 Hence, each function $u \in V(\o, q,\ell,\tau)$ can be written as
 \be u =  \a_0 \o +  \sum_{ 1 \leq i \leq q+\ell} \a_i \varphi_{a_i, \l_i}  + v \ee
 where the parameters $\a_i, a_i, \l_i$ are the solution of the previous minimization problem satisfying
  \be\label{alp6}
| \a_i^{4/(n-2)} K(a_i) - 1 | < c \, \tau \, \, \,  \forall \, \, 1 \leq i \leq q+\ell \quad  ; \quad  | \a_0 - 1 | < c\,  \tau \quad ; \quad \l_i d(a_i, \partial \mathbb{S}^n_+) > c \, \tau^{-1}  \, \, \, \forall i > q\ee
 and $v$ satisfies
 \be\label{wV0} \| v \| \leq \tau \, \,  ;   \quad  \langle v, \o \rangle = 0 \quad \mbox{  and } \quad  v \mbox{ satisfies \eqref{V0} }.\ee
 \end{pro}
Regarding the infinite dimensional part $v$ we prove the following proposition:
  \begin{pro} \label{wpv}
   Assume that $\o$ is a non degenerate solution of  $(\mathcal{P})$ and let $u:= \a_0 \o + \sum_{i=1}^{q+\ell} \a_i \varphi_{a_i,\l_i}  \in V(\o, q,\ell, \tau)$ with $\e \ln \l_i$ is small for each $i$. Then there exists a unique function $\ov{v}:= \ov{v}_{\e, \a_i, a_i,\l_i}$ such that
 \be\label{AAAZZZw}  \langle \n I_\e (u + \ov{v}) , v \rangle = 0 \quad  \mbox{ for each } v \mbox{ satisfying  \eqref{wV0}}.\ee
   Furthermore, it holds
$$ \|\overline{v} \| \leq c \, R(\e,a,\l) + c \,  R(\l) \quad \mbox{ with } R(\l) := \sum_{i=1}^{q+p}   \chi (\l_i) \quad\mbox{ where  } \chi(\l ) :=  \begin{cases}
 {1}/ {\l ^{3/2}} &  \mbox{ if } n = 5, \\
( {\ln^{2/3} \l } ) / {\l ^2} & \mbox{ if } n = 6,\\
 {1} / {\l ^2} & \mbox{ if } n \geq 7 . \end{cases}$$
(We remark that the term $R(\l)$ does not appear if $\o =0$).
   \end{pro}
\begin{pf}
Observe that,  for $b \geq 0$, $c \in \R$ and $\g > 2$, it holds
$$ | b+c |^{\g} = b^{\g} + \g b^{\g-1} c + \frac{1}{2}\g(\g-1) b^{\g-2}c^2 + O\big( |c|^{\g} + b^{\g-3} |c|^3 (\mbox{if } \g > 3)\big).$$
Thus we derive that, for each $v$ satisfying \eqref{wV0},
\begin{align*} & I_\e(u+v) = I_\e(u) - f_\e(v) + (1/2) Q_\e(v) + o(\| v\| ^2), \quad \mbox{ where } \\
& f_\e(v) = \int K u^{p-\e} v \quad \mbox{ and } \quad Q_\e(v):= \| v\| ^2 - ((n+2)/(n-2)) \int K u^{p-1-\e}v^2. \end{align*}
Moreover, using \eqref{alp6}, easy computations imply that
$Q_\e(v) = \ov{Q}_0(v) + o(\| v\|^2)$  where $$ \ov{Q}_0(v):=  \| v\| ^2 - \frac{n+2}{n-2} \int K \o^{p-1}v^2 - \frac{n+2}{n-2} \sum_{i=1}^{q+\ell} \int  \d_i^{p-1}v^2.$$
Since we assumed that $\o$ is a non degenerate solution of $(\mathcal{P}_0)$, arguing as in \cite{Bahri-Invariant} (pages 354-355) we derive that $\ov{Q}_0$ is a non degenerate quadratic form and so $Q_\e$. Hence the existence of $\ov{v}$ satisfying the equality of the proposition follows. Concerning the estimate of $\| \ov{v} \|$, we get that
$$ \| \ov{v}\| \leq c \| f_\e \|.$$
Note that, for each $v$ satisfying \eqref{wV0}, it holds (we denote by $\ov{u} := \sum \a_i \varphi_i$)
$$ f_\e(v) = \int K \ov{u}^{p-\e} v + \int K \o ^{p-\e}v + \sum O\Big( \int \d_i^{p-1} \o | v | + \int \d_i \o^{p-1} | v | \Big) .$$
Now we need to estimate the previous terms. For the second one, since $v \perp \o$, it holds
\be\label{vwwv} \int K \o ^{p-\e}v = \int K \o ^{p}v + O (\e \| v \| ) =  O (\e \| v \| ) .\ee
For the third and the fourth ones, observe that
$$ \Big( \int_{B(a,1)} \d_{a,\l}^{8n/(n^2-4)} \Big)^{(n+2)/(2n)} +  \Big( \int_{B(a,1)} \d_{a,\l}^{2n/(n+2)} \Big)^{(n+2)/(2n)}\leq c  \chi( \l) .$$
$$ \Big( \int_{\R^n \setminus B(a,1)} \d_{a,\l}^{2n/(n-2)} \Big)^{2/n} \leq \frac{c}{\l ^2} \qquad ; \qquad   \Big( \int_{\R^n \setminus B(a,1)} \d_{a,\l}^{2n/(n-2)} \Big)^{(n-2)/(2n)} \leq \frac{c}{\l ^{(n-2)/2}}.$$
Hence, using the Holder's inequalities and, in $B(a_i,1)$, the fact that $\o$ is $L^\infty$ bounded, we obtain
\be\label{vwdelta}  \int \d_i^{p-1} \o | v | + \int \d_i \o^{p-1} | v | \leq c \| v \| \chi(\l_i) \ee
It remains the first one. (We remark that in the case $\o=0$, this integral is the unique term in the linear form $f_\e$). Recall that we have $\e \ln \l_i$ is small for each $i$. It holds that
 \be\label{777} \int_{\mathbb{S}_+^n} K \ov{u}^{\frac{n+2}{n-2}-\e} v = \sum \a_i^{\frac{n+2}{n-2}-\e} \int_{\mathbb{S}_+^n} K \varphi_i ^{\frac{n+2}{n-2}-\e} v + O \Big( \sum_ {i\neq j } \int_{\mathbb{S}_+^n} \sup(\varphi_j, \varphi_i) ^{\frac{4}{n-2}} \inf(\varphi_j, \varphi_i)  | v | \Big).\ee
Observe that, for $n\geq 6$, it follows that $4/(n-2) \leq 1$. Hence, using Holder's inequality and \eqref{za1}, we get
\begin{align*}
\int_{\mathbb{S}_+^n} \sup(\varphi_j, \varphi_i) ^{\frac{4}{n-2}} \inf(\varphi_j, \varphi_i)  | v | & \leq   \int_{\mathbb{S}_+^n} (\varphi_j \varphi_i) ^{\frac{n+2}{2(n-2)}} | v | \leq c \| v \| \Big(  \int_{ \mathbb{S}_+^n} (\d_j \d_i) ^{\frac{n}{n-2}} \Big)^{\frac{n+2}{2n}}\\
& \leq c  \| v \| \e _{ij }^{\frac {n+2}{2(n-2)}}(\ln \e _{ij}^{-1})^{\frac{n+2}{2n}} \quad \mbox{ if } n \geq 6,\\
\int_{\mathbb{S}_+^5} \sup(\varphi_j, \varphi_i) ^{{4}/{3}} \inf(\varphi_j, P\d_i)  | v | & \leq c \| v \| \e _{ij }(\ln \e _{ij}^{-1})^{{3}/{5}} \mbox{ if } n = 5.
 \end{align*}
It remains to estimate the first integral in \eqref{777}. Recall that $v$ satisfies \eqref{wV0} which implies that $\int \d_i^p v =0$ for each $i \leq q+\ell$. Thus, it follows that
\begin{align} \Big| \int_{\mathbb{S}_+^n} K & \d_i ^{\frac{n+2}{n-2}-\e}  v \Big|  =  \Big| K(a_i) \int_{\mathbb{S}_+^n}  \d_i ^{\frac{n+2}{n-2}- \e} v +
O\Big( \int | K(x) - K(a_i) | \d_i ^p | v |  \Big) \Big|  \nonumber \\
& \leq c  \int | \d_i ^{-\e} - c_0^{-\e}\l_i^{-\e (n-2)/2} | \d_i ^p | v |  + c | \n K(a_i) | \int_ {\mathbb{S}^n_+} d( x,a_i)   \d_i ^{\frac{n+2}{n-2}} | v | +   c\int_ {\mathbb{S}^n_+} d( x,a_i)^2  \d_i ^{\frac{n+2}{n-2}} | v | \nonumber \\
& \leq c \, \| v \| \Big(\e +  \frac{| \n K(a_i) | }{\l _i} + \frac{1}{\l _i^2} \Big), \label{77723}\end{align}
where we have used the Holder's inequality and Lemma \ref{lowerL2}.\\
Thus the estimate of the first integral of \eqref{777} follows for $i \leq q$ (since in this case we have $\varphi_i=\d_i$).\\
 Now, for $i \geq q+1$ (that is $a_i \in \mathbb{S}^n_+$), we get
\be\label{23777}  \int_{\mathbb{S}_+^n} K \varphi_i ^{\frac{n+2}{n-2}-\e} v  =   \int_{\mathbb{S}_+^n} K \d_i ^{\frac{n+2}{n-2}-\e} v +  O\Big(\int_{\mathbb{S}_+^n} \d_i ^{\frac{4}{n-2}} | \varphi_i - \d_i | | v | \Big)\ee
The first integral is computed in \eqref{77723}. Finally, using the fact that  $| \varphi_i - \d_i | \leq c \min (1/(\l_i^{(n-2)/2} d_i^{n-2})\, ; \, \d_i)$ (see  Lemma \ref{lem:varphi}) and denoting by $B_i:= B(a_i, d_i)$, we get that
$$ \int_{\mathbb{R}_+^n \setminus B_i } \d_i ^{\frac{4}{n-2}} | \varphi_i - \d_i | | v | \leq c \int_{\mathbb{R}_+^n \setminus B_i } \d_i ^{\frac{n+2}{n-2}}  | v | \leq c \frac{ \| v \| } {(\l_id_i)^{(n+2)/2}}, $$
$$ \int_{ B_i } \d_i ^{\frac{4}{n-2}} | \varphi_i - \d_i | | v | \leq  c \begin{cases}
\frac{ 1 }{(\l_i d_i^2)^{3/2}} \int_{ B_i } \d_i ^{\frac{4}{3}}  | v | \leq c \frac{ \| v \| } {(\l_id_i)^{3}} \quad \mbox{ if } n = 5,\\
\frac{ 1 }{(\l_i d_i^2)^{(n+2)/4}} \int_{ B_i } \d_i ^{\frac{n+2}{2(n-2)}}  | v | \leq c  \| v \| \frac{ \ln^{(n+2)/(2n)} (\l_i d_i)}{(\l_i d_i)^{(n+2)/2}}  \quad \mbox{ if } n \geq 6. \end{cases} $$
The result follows.
\end{pf}

\section{Asymptotic expansion of the gradient}

\subsection{Expansion in the neighborhood at infinity $V(\o, q,\ell,\tau)$}

 In the following propositions, we will give the asymptotic expansion of the gradient in the set $V(\o, q,\ell,\tau)$. Note that, in this paper we do not need some precisions given in these propositions, but we present some details for possible future use.

\begin{pro}\label{devIeps3}
Let   $u:= \a_0 \o +\sum_{i=1}^{N}\a_i\varphi_{a_i,\l_i}  + \ov{v}\in V(\o, q,\ell,\tau) $ (with $N:= q+\ell$). Assume that $\e \ln \l_i$ is small for each $i$. For each $i \leq q$, there hold
\begin{align*}
&  \langle \n  I_\e(u),\d_i\rangle  =   \,  \a_iS_n \, \Big(1-\frac{\a_i^{p-1}K(a_i)}{ \l_i ^{\e(n-2)/2}} \Big) + O\Big( R_{1,i}(\e,a,\l)\Big) , \\
&  \langle \n I_\e(u),\l_i\frac{\partial\d_i}{\partial \l_i}\rangle  =  \frac{\a_i^{p}}{\l_i ^{\e(n-2)/2}} \Big( c_4 K(a_i) \, \e - \frac{c_3}{\l_i}\frac{\partial K}{\partial \nu}(a_i) \Big) + \frac{c_2}{2} \sum_{i\neq j \leq q}  {\a_j}{\l_i}\frac{\partial \e_{ij}}{\partial \l_i}  \Big( 1 - \frac{\a_i^{p-1}K(a_i)}{ \l_i ^{\e(n-2)/2}} \\
&   \qquad\qquad\qquad\quad   - \frac{\a_j^{p-1}K(a_j)}{ \l_j ^{\e(n-2)/2}}  \Big)
   - \ov{c}_6 \frac{\o (a_i) }{\l_i ^{(n-2)/2}} \a_0 \Big( 1 - \a_0^{p-1} - \frac{\a_i^{p-1}K(a_i)}{ \l_i^{\e (n-2)/2}} \Big)  +O\Big( \frac{1}{\l_i ^2 }+R_{2,i}(\e,a,\l)\Big) , \\
&   \langle \n I_\e(u),\frac{1}{\l_i}\frac{\partial\d_i}{\partial a_i}\rangle _{\lfloor T_{a_i}(\partial \mathbb{S}^n_+)}  =   \frac{c_2}{ 2 } \sum_{j \leq q; j \neq i } \frac{\a_j}{\l_i}\frac{\partial \e_{ij}}{\partial a_i}  \Big( 1 - \frac{\a_i^{p-1}K(a_i)}{ \l_i ^{\e(n-2)/2}} - \frac{\a_j^{p-1}K(a_j)}{ \l_j ^{\e(n-2)/2}}  \Big) \\
&   \qquad\qquad\qquad \qquad \qquad \quad - \frac{\a_i^{p}} {\l_i ^{\e (n-2)/2}} c_5 \frac{\n_T K(a_i)}{\l_i} + O\Big(  \frac{1}{\l_i ^2 } + R_{3,i}(\e,a,\l)\Big) ,
\end{align*}
where
\begin{align*} & R_{1,i}(\e,a,\l) := \e +  \frac{ | \n K(a_i) | }{\l_i } +  \frac{ 1}{\l_i^2 }   + \sum \e_{ij}  +  R^2(\e,a,\l) + O_{\o} \Big(  \frac{ 1 }{ \l_i ^{(n-2)/2}} +   \sum \frac{ \ln \l_k}{ \l_k ^{n/2}}  \Big) , \\
& R_{2,i}(\e,a,\l) :=  R_4(\e,a,\l) + \sum_{j \geq q+1}\e_{ij}  ,\\
& R_{3,i}(\e,a,\l) :=  R_4(\e,a,\l) + \sum_{j \geq q+1}\e_{ij}  + \sum_{j \leq q; j \neq i}   \l_j d (a_j , a_i ) \e_{ij}^{\frac{n+1}{n-2}}  , \\
 &   R_4(\e,a,\l) :=      R^2 (\e,a,\l) + \sum_{j\neq k}\e_{jk}^{\frac{n}{n-2}}\ln(\e_{jk}^{-1}) + O_{\o} \Big(  \sum \frac{ \ln \l_k}{ \l_k ^{n/2}} \Big)  ,\end{align*}
$$S_n := \int_{\R_+^n} \frac{ c_0^{2n/(n-2)} } { (1+|x|^2)^n}d x \, \, ; \, \, c_3:=  (n-2) c_0^{2n/(n-2)} \int_{\R^n_+} \frac{x_n (| x | ^2 - 1 )} { (1+|x|^2)^{n+1}}d x > 0 $$
$$ c_4:=   c_0^{2n/(n-2)} \frac{(n-2)^2}{8} \int_{\R^n} \frac{ (| x | ^2 - 1 )\ln(1+| x |^2)} { (1+|x|^2)^{n+1}}d x > 0 \, \, ;
\, \, c_2 := \int_{\R^n} \frac{ c_0^{2n/(n-2)} } { (1+|x|^2)^{(n+2)/2}}d x $$
$$ c_5:=  c_0^{2n/(n-2)} \frac{n-2}{n} \int_{\R^n} \frac{ | x | ^2 } { (1+|x|^2)^{n+1}}d x \, \, ; \, \,
 \ov{c}_6 :=  2^{(n-2)/2}   \frac{n-2}{4} \frac{ {c}_2 }{ c_0}$$ and where $R(\e,a,\l)$ is defined in Proposition \ref{pv},  $T_{a_i}(\partial \mathbb{S}^n_+)$ denotes the tangent space at the point $a_i$. \\
 (We remark that $O_{\o}(.)$ does not appear if $\o=0$.)
\end{pro}

\begin{pf}
From the definition of $I_\e$, it is easy to deduce that
\be \label{nablaI}\langle \n I_\e (u) , h \rangle  = \langle u , h \rangle  - \int_{\mathbb{S}^n_+ } K | u | ^{ p-1-\e} u h \qquad \mbox{ for each } h \in H^1(\mathbb{S}^n_+ ).\ee
We will focus on  the second assertion and we will give the principal changes for the other ones. Taking $h= \l_i \partial \d_i / \partial \l_i$, it holds: $ \langle \ov{v} , \l_i \partial \d_i / \partial \l_i \rangle = 0$. To estimate the first term of \eqref{nablaI}, using \eqref{91'} for $j \geq q+1$ and \eqref{91''} for $j \leq q$, there hold
\begin{align} \langle \varphi_j , \l_i \frac{\partial \d_i }{ \partial \l_i}  \rangle & =  O (\int  \d_j ^p \d_i) = O( \e_{ij} ) \quad \forall \, \, j\geq q+1 \quad , \quad \langle \d_i , \l_i \frac{\partial \d_i }{ \partial \l_i}  \rangle = \int_{ \R^n_+} \d_i ^p \l_i \frac{\partial \d_i }{ \partial \l_i}  = 0 , \label{po1}  \\
 \langle \d_j , \l _i \frac{\partial\d _i }{\partial\l _i} \rangle & =  \int_{\R^n_+} \d_j^\frac{n+2}{n-2}  \l _i \frac{\partial\d _i }{\partial\l _i} = \frac{1}{2}\int_{\R^n} \d_j^\frac{n+2}{n-2}  \l _i \frac{\partial\d _i }{\partial\l _i} = \frac{c_2}{2} \l_i \frac{\partial \e_{ij}}{\partial \l_i} + O\big(  \e _{ij }^{\frac{n}{n-2}}\ln (\e _{ij }^{-1})\big)  \mbox{ if } j \leq q, \label{po2}
 \end{align} 
 \begin{align} 
 \langle \o , \l_i \frac{\partial \d_i }{ \partial \l_i}  \rangle & =  p \int_{\mathbb{S}^n_+} \d_i^{p-1} \l_i \frac{\partial \d_i }{ \partial \l_i} \o = p \int_{\R^n_+} \d_i^{p-1} \l_i \frac{\partial \d_i }{ \partial \l_i} (\i\o) \nonumber \\
 & = (\i \o)({a}_i)\,  p \int_{\R^n_+} \d_i^{p-1} \l_i \frac{\partial \d_i }{ \partial \l_i} + O\Big( \int_{ \R^n_+} \d_i ^p | x- {a}_i| \Big),\nonumber\\
& = -  2^{(n-2)/2} \o(a_i) \frac{1}{2 c_0} \frac{n-2}{2} \frac{ {c}_2 }{\l_i ^{(n-2)/2}} + O\Big( \frac{1}{\l_i ^{n/2}} \Big) \quad  \mbox{ (using \eqref{4***}) } . \label{aaz4} \end{align}
where $(\i \o)$ is defined in \eqref{i44}. This completes the estimate of the first term of \eqref{nablaI} and we get
\be \label{pppp1}   \langle u , \l_i \partial \d_i / \partial \l_i \rangle = \frac{c_2}{2} \sum_{i \neq j \leq q} \a_j \l_i \frac{\partial \e_{ij}}{\partial \l_i}  - \ov{c}_6 \frac{\o(a_i) }{\l_i ^{(n-2)/2}} + O\Big( \frac{1}{\l_i ^{n/2}} +  \sum_{i\neq j \leq q}  \e _{ij }^{\frac{n}{n-2}}\ln (\e _{ij }^{-1}) + \sum_{j > q} \e_{ij}  \Big) .
\ee

Concerning the second part,  let $\ov{u}:= u-\ov{v}$, we have
\be\label{az5}  \int_{\mathbb{S}^n_+ } K | u | ^{ p-1-\e} u  \l_i \frac{\partial \d_i }{ \partial \l_i} =  \int_{\mathbb{S}^n_+ } K \ov{u}  ^{ p-\e}   \l_i \frac{\partial \d_i }{ \partial \l_i} + (p-\e) \int_{\mathbb{S}^n_+ } K \ov{u}  ^{ p-1-\e} \ov{v}  \l_i \frac{\partial \d_i }{ \partial \l_i} + O(\| \ov{v} \|^2) .\ee
For the second integral of \eqref{az5}, let $\O_i:= \{ x: \sum_{j\neq i} \a_j \varphi_j + \a_0 \o \leq \a_i \d_i /2\}$,  it holds
\begin{align}
 \int_{\mathbb{S}^n_+ } K \ov{u}  ^{ p-1-\e} \ov{v}  \l_i \frac{\partial \d_i }{ \partial \l_i}  = &  \int_{\mathbb{S}^n_+ } K (\a_i\d_i)  ^{ p-1-\e} \ov{v}  \l_i \frac{ \partial \d_i }{ \partial \l_i} + O\Big(\sum_{j\neq i} \int_{\O_i }\d_i^{p-1} \d_j | \ov{v} | + \int_{\mathbb{S}^n_+\setminus \O_i }\d_j^{p-1} \d_i | \ov{v} |\Big)  \label{az6}\\
 & + O\Big(   \int \d_i^{p-1} \o | \ov{v} | + \int \d_i \o^{p-1} | \ov{v} | \Big).\nonumber
   \end{align}
   The estimate of the last remainder term is given in \eqref{vwdelta}. \\
Observe that, if $n\geq 6$, then $p-1 \leq 1$ and therefore, using the Holder's inequality and \eqref{za1}, it holds
$$\sum_{j\neq i} \int_{\O_i }\d_i^{p-1} \d_j | \ov{v} | + \int_{\mathbb{S}^n_+\setminus \O_i }\d_j^{p-1} \d_i | \ov{v} |\leq c \sum_{j\neq i} \begin{cases}
&   \| \ov{v}  \| \, \e_{ij} \ln(\e_{ij}^{-1})^{3/5} \mbox{ if } n=5, \\
& \int (\d_i \d_j)^{p/2} | \ov{v} | \leq c  \| \ov{v} \| \e_{ij}^{p/2} \ln(\e_{ij}^{-1} )^{\frac{n+2}{2n}}\mbox{ if } n\geq 6 \end{cases} .$$

 Concerning the first integral of the right hand side of \eqref{az6}, we have
 $$  \int_{\mathbb{S}^n_+ } K \d_i ^{p-1-\e} \ov{v} \l_i \frac{\partial \d_i }{ \partial \l_i}  = \frac{c_0^{-\e}}{\l_i ^{\e(n-2)/2}}  \int_{\mathbb{S}^n_+ } K \d_i ^{p-1} \ov{v} \l_i \frac{\partial \d_i }{ \partial \l_i} + O\Big( \int_{\mathbb{S}^n_+ } | \d_i ^{-\e} - \frac{c_0^{-\e}}{\l_i ^{\e(n-2)/2}} |  \d_i ^{p} | \ov{v} |\Big) . $$
 For the first integral, since $\ov{v}$ satisfies \eqref{wV0} we need to  expand $K$ around $a_i$ and for the second one we will use the Holder's inequality and  Lemma \ref{lowerL2}. Thus we obtain

\begin{align}  \int_{\mathbb{S}^n_+ } K \d_i ^{p-1-\e} \ov{v} \l_i \frac{\partial \d_i }{ \partial \l_i} & = O\Big( | \n K(a_i) | \int_{\mathbb{S}^n_+}  d( x,a_i ) \d_i^p | \ov{v} | + \int_{\mathbb{S}^n_+}  d^2( x,a_i ) \d_i^p | \ov{v} | + \e\,  \| \ov{v} \|   \Big)\label{az7} \\
& = O\Big(  \| \ov{v} \| \Big( \e  + \frac{| \n K(a_i) |}{\l_i} + \frac{1}{\l_i ^2}\Big) \Big) \nonumber \end{align}
which completes the estimate of \eqref{az6} and we get (by combining the previous estimates) that
\be\label{fgh1} \int_{\mathbb{S}^n_+ } K \ov{u}  ^{ p-1-\e} \ov{v}  \l_i \frac{\partial \d_i }{ \partial \l_i}  = O(R^2(\e,a,\l) )+ O_\o( R^2(\l) )
\ee where $R(\e,a,\l)$ and   $R(\l)$ are defined in Propositions \ref{pv} and \ref{wpv} respectively. \\
It remains to estimate the first integral of the right hand side of \eqref{az5}.  Observe that, for $ 1 < \g \leq 3$ and $t_1, \cdots, t_{N+1} > 0$, it holds
\be\label{est11}
\Big| (\sum t_i)^\g - \sum t_i^\g - \g t_1^{\g-1} (\sum_{j\neq 1} t_j) \Big| t_1 \leq c \sum_{k\neq j } (t_k t_j)^{(\g + 1)/2} .\ee
Thus, for $n \geq 5$ (since $ 1 < p-\e := (n+2)/(n-2) - \e < 3$), using \eqref{est11}, it holds
\begin{align}
&  \int_{\mathbb{S}^n_+ }  K (\ov{u}) ^{p-\e} \l_i \frac{\partial \d_i }{ \partial \l_i}  = \sum_{j=1}^N  \int_{\mathbb{S}^n_+ } K (\a_j \varphi_j) ^{p-\e} \l_i \frac{\partial \d_i }{ \partial \l_i} +  \int_{\mathbb{S}^n_+ } K (\a_0 \o) ^{p-\e} \l_i \frac{\partial \d_i }{ \partial \l_i}\label{wdw2} \\
 & + (p-\e)  \int_{\mathbb{S}^n_+ } K (\a_i \d_i) ^{p-\e-1}(\sum_{j \neq i} \a_j \varphi_j + \a_0 \o) \l_i \frac{\partial \d_i }{ \partial \l_i} + \sum_{k\neq j}O\Big(\int (\d_k \d_j)^{n/(n-2)} + \int (\o \d_k)^{n/(n-2)} \Big).\nonumber
 \end{align}

Concerning the remainder term containing $\o$, it holds
\be\label{deltaw} \int_{\mathbb{S}^n_+  }  \d_{a,\l}^{\frac{n}{n-2}} \o^{\frac{n}{n-2}}  = \int_{\mathbb{R}^n_+  }  \d_{a,\l}^{\frac{n}{n-2}} (\i\o)^{\frac{n}{n-2}}  \leq c \int_{B(a,1)  }  \d_{a,\l}^{\frac{n}{n-2}}  + c \Big( \int _{\R^n\setminus B(a,1)} \d_{a,\l}^{2n/(n-2)} \Big) ^{1/2} \leq c \, \frac{\ln\l }{\l^{n/2}}  \ee
where $\i\o$ is defined in \eqref{i44} and the estimate of the other remainder term is given in \eqref{za1}.\\
In addition, for $j \geq q+1$, using \eqref{91'} and \eqref{lower}, we get
\be\label{dd2} |  \int_{\mathbb{S}^n_+ } K (\a_j \varphi_j) ^{p-\e} \l_i \frac{\partial \d_i }{ \partial \l_i} |  +  | \int_{\mathbb{S}^n_+ } K (\a_i \d_i) ^{p-\e-1} \a_j \varphi_j  \l_i \frac{\partial \d_i }{ \partial \l_i} | \leq c  \int ( \d_i ^p \d_j + \d_i \d_j ^p ) \leq c \e_{ij}.\ee
For $i\neq j \leq q$, we need to be more precise. Using \eqref{po2}, \eqref{za1}, Lemma \ref{lowerL2} and Holder's inequality we get

\begin{align}  \int_{\mathbb{S}^n_+ } K  \d_j ^{p-\e} \l_i \frac{\partial \d_i }{ \partial \l_i}   & = \frac{c_0^{-\e}}{ \l_j^{\e(n-2)/2}} \int_{\mathbb{S}^n_+} {K} \d_j^{p}  \l _i \frac{\partial\d _i }{\partial\l _i} +  O\Big(  \int_{\mathbb{S}^n_+} \d_j^{p-1} \Big| \d_j ^{-\e} - \frac{c_0^{-\e}}{ \l_j^{\e(n-2)/2}}\Big| (\d_j \d_i) \Big)  \label{wxcv1} \\
&=\frac{ K(a_j) }{ \l_j^{\e(n-2)/2}}  \int_{ B_j }  \d_j^{p}  \l _i \frac{\partial\d _i }{\partial\l _i} + O\Big(  | \n \tilde{K}(a_j) |  \int_{ B_j  }| x-a_j|  \d_j^{p-1}  (\d_j \d_i)  \nonumber\\
 &  + \int_{ B_j  }  |x-a_j|^2  \d_j^{p-1}  (\d_j \d_i) + \int_{ \R^n_+\setminus B_j } \d_j ^{p-1} (\d_j \d_i) +  \e  \e_{ij}(\ln \e _{ij }^{-1}  )^{\frac{n-2}{n}}\Big)  \nonumber \\
 & = \frac{ K(a_j) }{ \l_j^{\e (n-2)/{2}}}   \langle \d_j ,  \l _i \frac{\partial\d _i }{\partial\l _i} \rangle +  O\Big(  \Big(\frac{| \n K(a_j) | }{\l_j}  + \frac{\ln \l_j  }{\l_j^2} +  \e \Big) \e_{ij}(\ln \e _{ij }^{-1}  )^{\frac{n-2}{n}}\Big)  \nonumber
 \end{align}
 where $B_j := B(a_j,1) \cap \R^n_+$.
 Following the same computations, for $i \neq j \leq q$, we get
\be p  \int_{\mathbb{S}^n_+ } K  \d_j \d_i ^{p-\e-1} \l_i \frac{\partial \d_i }{ \partial \l_i}    =   \frac{ K(a_i) }{ \l_i^{\e (n-2)/{2}}}   \langle \d_j ,  \l _i \frac{\partial\d _i }{\partial\l _i} \rangle +  O\Big(  \Big(\frac{| \n K(a_i) | }{\l_i}  + \frac{\ln \l_i  }{\l_i^2} +  \e \Big) \e_{ij}(\ln \e _{ij }^{-1}  )^{\frac{n-2}{n}}\Big)  .\ee

 Concerning the first term of \eqref{wdw2} for $k =i$,  using Lemma \ref{lowerL2}, we get
\be \label{dd1} \int_{\mathbb{S}^n_+ } K  \d_i ^{ p-\e}   \l_i \frac{\partial \d_i }{ \partial \l_i}  =\frac{ c_0^{-\e} }{\l_i^{\e\frac{n-2}{2}}} \Big( \int_{\mathbb{S}^n_+ } K  \d_i  ^{ p}   \l_i \frac{\partial \d_i }{ \partial \l_i} + \e \frac{n-2}{2}  \int_{\mathbb{S}^n_+ } K  \d_i  ^{ p}   \l_i \frac{\partial \d_i }{ \partial \l_i}\ln[2+ (\l_i^2-1)(1-\cos d(  x,a_i ))] \Big) + O( \e^2 ).
\ee
In addition, denoting $\tilde{K}:= K\circ \pi_{-a_i}$ (defined in \eqref{chvar}) and $(x-{a}_i)_n$ the $n$-th component of the vector $x-{a}_i \in \R^{n}$. Using \eqref{po1} and the change of variables defined in \eqref{chvar} ,  there hold
\begin{align}
 \int_{\mathbb{S}^n_+ } K  \d_i  ^{ p}   \l_i \frac{\partial \d_i }{ \partial \l_i}  & = \int_{\R^n_+ } \tilde{K}  \d_i  ^{ p}   \l_i \frac{\partial \d_i }{ \partial \l_i} =  \frac{\partial \tilde{K}}{\partial x_n}( {a}_i) \int_{\R^n_+ } (x-{a}_i)_n  \d_i  ^{ p}   \l_i \frac{\partial \d_i }{ \partial \l_i} + O\Big( \int_{\R^n_+ } | x-{a}_i |^2  \d_i  ^{ p+1} \Big) \nonumber \\
&  = - 2 \frac{\partial K}{\partial \nu} (a_i) \frac{n-2}{2} \frac{c_0^{2n/(n-2)}}{\l_i}\int_{\R^n_+} \frac{ y_n (1-| y |^2)}{(1+ | y |^2)^{n+1}}dy + O(\frac{1}{\l_i^2}), \label{ppp3} \end{align}
\begin{align}
& \int_{\mathbb{S}^n_+ } K  \d_i  ^{ p}   \l_i \frac{\partial \d_i }{ \partial \l_i}\ln [ 2+ (\l_i^2-1)(1-\cos d ( x,a_i )) ]  =   \int_{\mathbb{R}^n_+ } \tilde{K}  \d_{0,\l_i}  ^{ p}   \l_i \frac{\partial \d_{0,\l_i} }{ \partial \l_i}\ln [ \frac{1}{1+|y|^2}(2+2\l_i^2|y|^2) ]  \nonumber\\
& \qquad \qquad \qquad  = \int_{\mathbb{R}^n_+ } \tilde{K}  \d_{0,\l_i}  ^{ p}   \l_i \frac{\partial \d_{0,\l_i} }{ \partial \l_i}\ln (2+2\l_i^2|y|^2) + O \Big( \int_{\mathbb{R}^n_+ } \d_{0,\l_i}  ^{ p+1}  \ln (1+|y|^2) \Big) \label{rem1}\\
& \qquad \qquad \qquad  = K(a_i) c_0^\frac{2n}{n-2} \frac{n-2}{2}  \int_{\R^n_+} \frac{ (1-| y |^2)\ln(1+ | y |^2)}{(1+ | y |^2)^{n+1}}dy +O(\frac{1}{\l_i}), \label{ppp4}
\end{align}
which achieve  the estimate of \eqref{dd1} (note that the estimate of the remainder term in \eqref{rem1} is done in \eqref{app:1} by using the fact that $\ln(1+|y|^2) \leq |y|^2$).
It remains to estimate the integrals of \eqref{wdw2} which contain $\o$. Using \eqref{az1} and Lemma \ref{lowerL2}, there hold
\begin{align}
& \int_{\mathbb{S}^n_+ } K  ( \a_0 \o )  ^{ p-\e}   \l_i \frac{\partial \d_i }{ \partial \l_i} = \a_0^{p} \int K\o ^p  \l_i \frac{\partial \d_i }{ \partial \l_i} + O\Big( \e\, \int K \o^p \d_i\Big) = \a_0 ^p \langle \o ,  \l_i \frac{\partial \d_i }{ \partial \l_i} \rangle + O\Big( \frac{ \e }{ \l_i ^{(n-2)/2}} \Big), \label{fgh2}\\
& p \int_{\mathbb{S}^n_+ } K   \d_i ^{ p-\e-1}  \o   \l_i \frac{\partial \d_i }{ \partial \l_i} \label{fgh3} \\
& = \frac{c_0^{-\e}K(a_i) }{\l_i^{\e(n-2)/2}}  p \int_{\mathbb{S}^n_+ }    \d_i ^{ p-1}   \l_i \frac{\partial \d_i }{ \partial \l_i} \o + O\Big( \int_{\mathbb{S}^n_+ } d(x,a_i) \d_i^p  + \e \int_{\mathbb{S}^n_+ } \ln [2+(\l_i^2-1)(1-\cos d( x,a_i ))] \d_i^p  \Big) \nonumber \\
& = \frac{c_0^{-\e}K(a_i) }{\l_i^{\e(n-2)/2}} \,  \langle \o ,  \l_i \frac{\partial \d_i }{ \partial \l_i} \rangle + O\Big(\frac{ 1 }{ \l_i ^{n/2}}  + \frac{ \e }{ \l_i ^{(n-2)/2}} \Big). \nonumber
\end{align}

Hence, combining \eqref{za1}, \eqref{deltaw}-\eqref{ppp3} and \eqref{ppp4}-\eqref{fgh3}, we derive the estimate of \eqref{wdw2} and we get
\begin{align}\label{ppp5}
 \int_{\mathbb{S}^n_+ }  K (\ov{u}) ^{p-\e} \l_i \frac{\partial \d_i }{ \partial \l_i}  = &   \sum_{i\neq j \leq q}  \a_j \Big( \frac{ \a_i ^{p-1} K(a_i) }{ \l_i^{\e (n-2)/{2}}}  + \frac{ \a_j ^{p-1} K(a_j) }{ \l_j^{\e (n-2)/{2}}}\Big)  \langle \d_j ,  \l _i \frac{\partial\d _i }{\partial\l _i} \rangle + c_3 \frac{\a_i p }{\l_i^{\e(n-2)/2}} \frac{\partial K}{\partial \nu} (a_i) \\
 & + \a_0 \Big(  \a_0^{p-1} -  \frac{  \a_i ^{p-1}K(a_i) }{\l_i^{\e(n-2)/2}} \Big) \frac{ \ov{c}_6 \, \o(a_i) }{\l_i ^{(n-2)/2}} + O\Big( \frac{1}{\l_i ^2 }+R_{2,i}(\e,a,\l)\Big).
\nonumber \end{align}
Finally, combining \eqref{nablaI}, \eqref{po2}, \eqref{pppp1}, \eqref{az5}, \eqref{az7}, \eqref{ppp5} and the estimate of $\| \ov{v}\|$ (given in Proposition \ref{wpv}), the second claim follows.

To prove the first claim, in \eqref{nablaI}, we consider $h=\d_i$. Observe that, since $\ov{v}$ satisfies \eqref{wV0}, we get $ \langle \ov{v} , \d_i \rangle = 0$. Furthermore, using \eqref{i44}-\eqref{eq:delta}, \eqref{91'}, \eqref{3***} and the fact that $\o$ is $L^\infty$ bounded, we obtain
\begin{align}
& \langle \varphi_j , \d_i \rangle = \int_{\R^n_+} \d_j ^p \d_i = O( \e_{ij} ) \quad \forall \, \,  j\neq i,  \quad  \langle \o , \d_i \rangle = \int_{\R^n_+} \d_i ^p (\i\o) = O( \frac{1}{\l_i^{(n-2)/2}}), \label{az1}\\
& \| \d_i \| ^2 = \int_{ \mathbb{S}^n_+ } \d_i ^{p+1} = \int_{\R^n_+} \d_i ^{p+1} = c_0^{p+1} \int_{\R^n_+} \frac{ dx}{ (1+| x |^2)^n} := S_n \label{az2}
\end{align}
which give the first part of \eqref{nablaI}.
For the second one, we note that \eqref{az5} and \eqref{fgh1} hold with $\d_i$ instead of $\l_i \partial \d_i /\partial \l_i$. Moreover, since the $\e \ln \l_k $'s are small, it holds
$$ \int_{\mathbb{S}^n_+ } K \ov{u}^{p -\e}  \d_i = \int_{\mathbb{S}^n_+ } K \d_i ^{p+1-\e} + O\Big( \int \d_i^{p} (\sum_{j \neq i} \d_j + \o) + \int \d_i ( \sum_{j \neq i } \d_j ^p + K \o ^p) \Big).$$
Using \eqref{az1} and \eqref{91'}, the remainder term is dominated by $ \sum_{j \neq i} \e_{ij} + 1/\l_i ^{(n-2)/2}$. For the other integral,
since we assumed that $\e \ln \l_i $ is small then using Lemma  \ref{lowerL2}  and expanding ${K}$ around $a_i$, we get
\begin{align} \int_{\mathbb{S}^n_+ } K \d_i ^{p+1-\e} & =  K(a_i)  \frac{ c_0^{-\e} }{ \l_i^{\e\frac{n-2}{2}}} \int_{ \mathbb{S}^n_+ }  \d_i ^{p+1}  + O \Big(\int_{\mathbb{S}^n_+} | K( x) - K(a_i)  | \d_i ^{p+1}  + \int _{\mathbb{S}^n_+ } | \d_i^{-\e} - \frac{ c_0^{-\e} }{ \l_i^{\e\frac{n-2}{2}}}  |  \d_i ^{p+1} \Big)\label{zza1}\\
& = K(a_i) \l_i^{-\e(n-2)/2} S_n + O\Big( \e + \frac{ | \n K(a_i) | }{\l_i } +  \frac{1}{\l_i ^2 }  \Big)\nonumber\end{align}
by  using \eqref{az2} and the fact that $c_0^{-\e}  = 1 + O(\e)$. Thus the first claim follows.

Concerning the last claim, it follows as the second one by taking $h= (1/\l_i) \partial \d_i / \partial a_i$.  We recall that the concentration point $a_i$ belongs to $\partial \mathbb{S}^n_+$. Hence its image via the stereographic projection will be on the boundary of $\R^n_+$ that is its $n$-th  component is $0$. Moving the point on the boundary will imply $(n-1)$ derivatives with respect to $a_i^k$ with $1 \leq k \leq n-1$.

Observe that, using the second assertion of Estimate $F11$ (page 22) of \cite{B1}, there hold
\begin{align}
 \langle \d_j , \frac{1}{\l_i} \frac{\partial \d_i }{ \partial a_i}  \rangle & = \int_{\R^n_+}  \d_j ^p \frac{1}{\l_i} \frac{\partial \d_i }{ \partial a_i} = \frac{1}{2} \int_{\R^n}  \d_j ^p \frac{1}{\l_i} \frac{\partial \d_i }{ \partial a_i} = \frac{1}{2}  c_2\frac{1}{\l_i} \frac{\partial \e_{ij} }{ \partial a_i} + O(\l_j d( a_i , a_j ) \e_{ij} ^{\frac{n+1}{n-2}})  \mbox{ if } a_j \in \partial \mathbb{S}^n_+,\nonumber \\
 \langle \varphi_j , \frac{1}{\l_i} \frac{\partial \d_i }{ \partial a_i}  \rangle & = \int_{\R^n_+}  \d_j ^p \frac{1}{\l_i} \frac{\partial \d_i }{ \partial a_i} = O \Big( \int _{\R^n}\d_j ^p \d_i \Big) = O(\e_{ij} ) \mbox{ if } a_j \in \mathbb{S}^n_+,\nonumber  \\
 \langle \d_i , \frac{1}{\l_i} \frac{\partial \d_i }{ \partial a_i}  \rangle & = \int_{\R^n_+}  \d_i ^p \frac{1}{\l_i} \frac{\partial \d_i }{ \partial a_i} = 0, \nonumber  \\
 \langle \o , \frac{1}{\l_i} \frac{\partial \d_i }{ \partial a_i}  \rangle & = p \int_{\mathbb{S}^n_+}  \d_i ^{p-1} \frac{1}{\l_i} \frac{\partial \d_i }{ \partial a_i} \o   =p \int_{\R^n_+}  \d_i ^{p-1} \frac{1}{\l_i} \frac{\partial \d_i }{ \partial a_i}( \i \o) = O \Big( \int \d_i ^p | x - a_i | \Big) = O\Big( \frac{1}{\l_i ^{n/2}} \Big),\label{www11} \end{align}
 which give the estimate of the first part of \eqref{nablaI}. Concerning the second one, \eqref{az5} holds with $(1/\l_i) \partial \d_i / \partial a_i$ instead of $\l_i \partial \d_i / \partial \l_i$. Furthermore, the estimate of the second integral of \eqref{az5} holds in the same way as before which implies that \eqref{fgh1} holds with $(1/\l_i) \partial \d_i / \partial a_i$ instead of $\l_i \partial \d_i / \partial \l_i$. Concerning the first integral of \eqref{az5}, we expand it as in \eqref{wdw2} by taking $(1/\l_i) \partial \d_i / \partial a_i$ instead of $\l_i \partial \d_i / \partial \l_i$. \\
$ (*) $  For the  terms containing $\o$  of the analogue equation of \eqref{wdw2}, using \eqref{lower}, \eqref{az1} and \eqref{www11}, it holds
\begin{align*}
& \int_{\mathbb{S}^n_+ }  K \o ^{p-\e}  \frac{1}{ \l_i} \frac{\partial \d_i }{ \partial a_i} = \int_{\mathbb{S}^n_+ }  K \o ^{p}  \frac{1}{ \l_i} \frac{\partial \d_i }{ \partial a_i } + O\Big( \e \int K \o^p \d_i \Big) =  \langle \o   \frac{1}{ \l_i} \frac{\partial \d_i }{ \partial a_i } \rangle + O \Big( \frac{\e}{\l_i^{(n-2)/2}} \Big) = O\Big( \e^2 + \frac{1}{\l_i ^{n/2}}\Big),\\
& p \int_{\mathbb{S}^n_+ } K  \d_i  ^{ p-\e-1}  \frac{1}{ \l_i} \frac{\partial \d_i }{ \partial a_i} \o  =   p \int_{\mathbb{S}^n_+ } K  \d_i  ^{ p-1}  \frac{1}{ \l_i} \frac{\partial \d_i }{ \partial a_i} \o + O \Big(\e \ln \l_i  \int \d_i ^p \o \Big) \\
& \qquad \qquad \qquad \qquad \qquad  = K(a_i) \langle \o,  \frac{1}{ \l_i} \frac{\partial \d_i }{ \partial a_i} \rangle + O\Big( \int_{\R^n_+}  d( x ,a_i ) \d_i^{p} + \e \frac{\ln \l_i}{\l_i ^{(n-2)/2}} \Big)  = O\Big( \frac{1}{\l_i ^{n/2}} + \e ^2 \Big).
\end{align*}

$ (*) $ Furthermore,  for $j > q$ that is $a_j \in \mathbb{S}^n_+$, since $\e \ln \l_j$ and  $\e \ln\l_i$ are small, it holds
$$ \Big| \int_{\mathbb{S}^n_+ }  K \varphi_j ^{ p-\e}  \frac{1}{ \l_i} \frac{\partial \d_i }{ \partial a_i^k} \Big| +
\Big| p \int_{\mathbb{S}^n_+ } K  \d_i  ^{ p-\e-1} \varphi_j \frac{1}{ \l_i} \frac{\partial \d_i }{ \partial a_i}\Big|  \leq c \Big( \int \d_i ^p \d_j +  \int \d_i  \d_j ^p \Big) = O(\e_{ij})$$
which gives the estimate of a part of the first and the third terms of the analogue equation of \eqref{wdw2}.\\
For the other terms, we need to be more precise. Taking $k \leq n-1$,

$ (*) $ For the first integral of the analogue equation of  \eqref{wdw2} with $j = i$,  using Lemma \ref{lowerL2}, it holds
\begin{align*}
 & \int_{\mathbb{S}^n_+ }  K \d_i  ^{ p-\e}  \frac{1}{ \l_i} \frac{\partial \d_i }{ \partial a_i^k} \\
& = \frac{ c_0^{-\e} }{\l_i^{\e(n-2)/2}} \Big( \int_{ \mathbb{S}^n_+  } {K}  \d_i  ^{ p}  \frac{1}{ \l_i} \frac{\partial \d_i }{ \partial a_i^k} + \e \frac{n-2}{2} \int_{ \mathbb{S}^n_+  } {K}  \d_i  ^{ p} \ln [ 2 +(\l_i^2-1)(1-\cos d ( x,a_i) )] \frac{1}{ \l_i} \frac{\partial \d_i }{ \partial a_i^k} \Big) +O ( \e ^2  ) .
\end{align*}
For the second integral, using the change of variable defined in \eqref{chvar}, the parity of the function ${\partial \d_i }/{ \partial a_i^k}$ for $k \leq n-1$ and expanding $\tilde K:= K \circ \pi_{-a_i}$ around $0$, we get
$$  \int_{ \mathbb{S}^n_+  } {K}  \d_i  ^{ p} \ln [ 2 +(\l_i^2-1)(1-\cos d ( x,a_i) )] \frac{1}{ \l_i} \frac{\partial \d_i }{ \partial a_i^k} = \int_{ \mathbb{R}^n_+  } \tilde{K}  \d_i  ^{ p} \ln [ \frac{2}{1+|y|^2}(1+ \l_i^2|y|^2) ] \frac{1}{ \l_i} \Big(\frac{\partial \d_i }{ \partial {a}_i^k}\Big)_{| {a}_i =0}  = O\Big( \frac{1}{\l_i} \Big).$$
Concerning the first one, using again the change of variable defined in \eqref{chvar}, it follows that (since $k \leq n-1$)
\begin{align*} \int_{ \mathbb{S}^n_+  } {K}  \d_i  ^{ p}  \frac{1}{ \l_i} \frac{\partial \d_i }{ \partial a_i^k} & =  \int_{ \mathbb{R}^n_+  } \tilde{K}  \d_{0,\l_i}  ^{ p}  \frac{1}{ \l_i} \Big(\frac{\partial \d_i}{\partial {a}_i^k}\Big)_{|{a}_i =0} = c_0^{p+1}(n-2) \int_{\mathbb{R}^n_+} \tilde{K} \frac{\l^{n+1} y_k}{(1+\l_i ^2 | y |^2)^{n+1}}dy \\
&  = \frac{\partial \tilde{K}}{\partial x^k}(0) c_0^{p+1}(n-2) \frac{1}{2n} \frac{1}{\l_i}\int_{\R^n } \frac{ | y |^2}{(1+| y |^2)^{n+1}} dy + O\big( \frac{1}{\l_i ^2}   \big).\end{align*}

$ (*) $  For  $i \neq j \leq q$ that is $a_j \in \partial \mathbb{S}^n_+ $,  using \eqref{za1}, Lemma \ref{lowerL2} and the Holder's inequality, following the computations done in \eqref{wxcv1}, there hold

\begin{align*}
& \int_{\mathbb{S}^n_+ } K   \d_j  ^{ p-\e}  \frac{1}{ \l_i} \frac{\partial \d_i }{ \partial a_i}   = \frac{1}{\l_j^{\e(n-2)/2}}K(a_j) \langle \d_j,  \frac{1}{ \l_i} \frac{\partial \d_i }{ \partial a_i} \rangle + O\Big( \e_{ij} | \ln\e_{ij} |^{\frac{n-2}{n}} (\frac{ | \n K(a_j) |  }{\l_j} + \frac{\ln\l_j}{\l_j^2} + \e ) \Big),\\
& p \int_{\mathbb{S}^n_+ } K  \d_i  ^{ p-\e-1} \d_j   \frac{1}{ \l_i} \frac{\partial \d_i }{ \partial a_i} = \frac{1}{\l_i^{\e(n-2)/2}}  K(a_i) \langle \d_j,  \frac{1}{ \l_i} \frac{\partial \d_i }{ \partial a_i} \rangle + O\Big( \e_{ij} | \ln\e_{ij} | ^{\frac{n-2}{n}}( \frac{ | \n K(a_i) | }{\l_i} + \frac{\ln \l_i}{\l_i ^2} + \e )  \Big).
\end{align*}

The proof of the third assertion follows. This ends the proof of the proposition.
\end{pf}

 \begin{pro}\label{devIeps4}
Let   $u:= \a_0 \o +\sum_{i=1}^{N}\a_i\varphi_{a_i,\l_i}  + \ov{v}\in V(\o, q,\ell,\tau) $ (with $N:= q+\ell$). Assume that  $\e \ln \l_i$ is small for each $i$. Then, for each $i > q$, there hold

\begin{align*}
& \langle \n  I_\e(u),\o \rangle =  \a_0\| \o \|^2 \, (1-\a_0^{p-1})+O\Big( \e + \| \ov{v} \|^2 + \sum \frac{1}{\l_i ^{(n-2)/2}}+ \frac{1}{\l_i ^4} \Big) ,\\
& \langle \n  I_\e(u),\varphi_i\rangle = 2 \a_iS_n \, \Big(1- \frac{\a_i^{p-1}K(a_i)}{ \l_i^{\e(n-2)/2}} \Big)+O\Big(\e  + \sum \e_{ij} +   \frac{ 1}{\l_i ^2 } +  R^2(\e,a,\l) + \frac{1}{(\l_i d_i)^{n-2}}  \Big)\\
& \qquad \qquad \qquad + O_\o\Big( \frac{1}{\l_i ^{(n-2)/2} } + \sum \frac{\ln \l_k}{\l_k ^{n/2}} \Big), \\
 & \langle \n I_\e(u),\l_i\frac{\partial\varphi_i}{\partial \l_i}\rangle = \frac{  \a_i^{p} }{\l_i ^{\e(n-2)/2}} \Big( 2 c_4 K(a_i) \, \e + {c_7}\frac{\D K(a_i)}{\l_i^2} \Big)  + {c_2} \sum_{j \neq i}  {\a_j}{\l_i}\frac{\partial \e_{ij}}{\partial \l_i}  \Big( 1 - \frac{\a_i^{p-1}K(a_i)}{ \l_i ^{\e(n-2)/2}} \\
 & \qquad \qquad \qquad - \frac{\a_j^{p-1}K(a_j)}{ \l_j ^{\e(n-2)/2}}  \Big) - {c_2}\frac{n-2}{2} \sum_{ j > q}  {\a_j} \frac{ H(a_j,a_i)}{(\l_i \l_j)^{n-2)/2}}  \Big( 1 - \frac{\a_i^{p-1}K(a_i)}{ \l_i ^{\e(n-2)/2}} - \frac{\a_j^{p-1}K(a_j)}{ \l_j ^{\e(n-2)/2}}  \Big) \\
 & \qquad \qquad \qquad - 2 \ov{c}_6 \frac{\o(a_i)}{\l_i^{(n-2)/2}} \a_0 \Big( 1 - \a_0^{p-1} - \frac{\a_i^{p-1}K(a_i)}{ \l_i^{\e (n-2)/2}} \Big)
 + O\Big(R_4(\e,a,\l)  + \sum_{ j > q} \frac{\ln(\l_k d_k)}{(\l_k d_k)^n} \Big)\\
 &  \langle \n I_\e(u),\frac{1}{\l_i}\frac{\partial \varphi_i}{\partial a_i}\rangle =  - 2\a_i^{p}c_5\frac{\n K(a_i)}{\l_i} + O\Big(  \frac{1}{\l_i ^3} + \frac{1}{(\l_i d_i)^{n-2}}+ \sum \e_{ik} + R^2(\e,a,\l) \Big) + O_\o\Big( \sum \frac{\ln \l_k}{\l_k ^{n/2}} \Big)
\end{align*}
where the constant $c_7$ is defined in \eqref{c7}, $R_4(\e,a,\l)$ and the other constants  are introduced in Proposition \ref{devIeps3} and $R(\e,a,\l)$ is defined in Proposition \ref{pv}.
\end{pro}

\begin{pf}
For the first one, since $\ov{v}$ satisfies \eqref{wV0}, using \eqref{nablaI} and \eqref{az1}, it holds
$$ \langle \n  I_\e(u),\o \rangle = \langle u,\o \rangle - \int K | u |^{p-\e-1} u \o = \a_0 \| \o \| ^2 + \sum O \Big( \frac{1}{\l_i ^{(n-2)/2}} \Big) - \int K | u |^{p-\e -1} u \o $$ and, as in the proof of \eqref{az5}, we have $$\int_{\mathbb{S}^n_+ } K | u | ^{ p-1-\e} u  \o =  \int_{\mathbb{S}^n_+ } K \ov{u}  ^{ p-\e}   \o + (p-\e) \int_{\mathbb{S}^n_+ } K \ov{u}  ^{ p-1-\e} \ov{v}  \o  + O(\| \ov{v} \|^2) .$$
Note that, as in the proof of \eqref{fgh1}, it follows that $$ \int_{\mathbb{S}^n_+ } K \ov{u}  ^{ p-1-\e} \ov{v}  \o = \a_0^{p-1-\e} \int_{\mathbb{S}^n_+ } K \o ^{p-\e} \ov{v} + \sum_{j\geq 1} O\Big( \int \o  ^{p-1} \d_j | \ov{v} | + \int \o \d_j ^{p-1} | \ov{v} | \Big). $$
Recall that these terms are computed in \eqref{vwwv} and  \eqref{vwdelta} respectively. Furthermore, using \eqref{az1}, it holds
$$ \int_{\mathbb{S}^n_+ } K \ov{u}  ^{ p-\e}   \o  = \int_{\mathbb{S}^n_+ } K (\a_0 \o)  ^{ p-\e}   \o +  \sum_{j\geq 1} O\Big( \int \o  ^{p} \d_j + \int \o \d_j ^{p}  \Big) =
\a_0   ^{ p} \| \o \| ^2 + O\Big( \e + \sum_{j\geq 1}  \frac{ 1 }{ \l_j ^{(n-2)/2} } \Big).$$
Hence the proof of the first claim follows.

Concerning the other claims, they follow as the previous proposition by taking $\varphi_i$ instead of $\d_i$ and we need to use Lemma \ref{lem:varphi}. We start by the following remark: for $a_i \in \mathbb{S}^n_+$ with $\l_i d_i$ large, we have
$$  \int_{\mathbb{S}^n_+ } \d_i ^{p+1} =  \int_{\mathbb{S}^n } \d_i ^{p+1} + O\Big(  \int_{\mathbb{S}^n \setminus B(a_i, d_i) } \d_i ^{p+1} \Big)= 2 S_n +  O\Big( \frac{1}{(\l_i d_i )^n} \Big) $$ which implies, comparing with Proposition \ref{devIeps3}, that some constants will be multiplied by 2.

Precisely,  for the second claim, \eqref{az1} holds with $\varphi_i$ instead of $ \d_i$ and  \eqref{az2} has to be changed as follows
$$ \| \varphi_i\|^2 = \int_{\mathbb{S}^n_+} \d_i ^p \varphi_i  = \int_{ \mathbb{S}^n_+ } \d_i ^{p+1} + O\Big( | \varphi_i - \d_i |_\infty \int_{\R^n} \d_i ^{p}  \Big) = 2 S_n + O\Big(\frac{1}{(\l_i d_i)^{n-2}} \Big)$$
(where we have used Lemma \ref{lem:varphi}) which give the estimate of the first part of \eqref{nablaI}. Concerning the second part, \eqref{az5} holds with $\varphi_i$ instead of $\l_i \partial  \d_i/ \partial \l_i$. The expansion of the analogue of \eqref{az6} will contain
$$ \int_{\mathbb{S}^n_+ } K (\a_i\varphi_i)  ^{ p-1-\e} \ov{v}  \varphi_i \quad \mbox{ instead of } \quad \int_{\mathbb{S}^n_+ } K (\a_i\d_i)  ^{ p-1-\e} \ov{v}  \l_i \frac{ \partial \d_i }{ \partial \l_i}$$ and this term is computed in \eqref{23777}. Furthermore, \eqref{zza1} will be changed as follows
$$ \int_{\mathbb{S}^n_+} K\varphi_i^{p+1-\e} = \int_{B(a_i, d_i)} K\d_i^{p+1-\e} + O\Big(\frac{1}{(\l_i d_i)^{n-2}} \Big) = 2 S_n K(a_i) + O\Big( \e + \frac{1}{\l_i ^2} + \frac{1}{(\l_i d_i)^{n-2}} \Big)$$
(by using Lemmas \ref{lowerL2} and \ref{lem:varphi}). Hence the second assertion follows.

Concerning the third assertion of the proposition, it follows as the corresponding one in Proposition \ref{devIeps3} with some changes since the point $a_i \in \mathbb{S}^n_+$ and therefore we will have $\varphi_i$ instead of $\d_i$. Using \eqref{za1} and \eqref{91''},  we get

\begin{align}
p  \int_{\mathbb{S}^n_+} \d_i ^{p-1}{\l _i} \frac{\partial\d_i }{\partial \l _i } \d_j & = p  \int_{\R^n_+} \d_i ^{p-1}{\l _i} \frac{\partial\d_i }{\partial \l _i } \d_j = p  \int_{\R^n} \d_i ^{p-1}{\l _i} \frac{\partial\d _i }{\partial \l _i } \d_j + O\Big( \int_{\R^n\setminus B(a_i,d_i)} \d_i ^p \d_j\Big) \label{qq11}  \\
& = c_2 \l_i\frac{\partial\e_{ij}}{\partial\l_i} + O\Big( \e _{ij }^{\frac{n}{n-2}}  \ln \e _{ij }^{-1} + \frac{1}{(\l_i d_i)^2} \e_{ij} (\ln \e _{ij }^{-1}  )^{\frac{n-2}{n}}\Big). \nonumber \end{align}
Hence the first term of \eqref{nablaI}, for $j \leq q$, follows. Now, for $j \geq q+1$ with $j\neq i$, using Lemma \ref{lem:varphi}, we get
\be\label{qq112} \langle \varphi_j,{\l _i} \frac{\partial\varphi _i }{\partial \l _i } \rangle = p  \int_{\R^n_+} \d_i ^{p-1}{\l _i} \frac{\partial\d_i }{\partial \l _i } \d_j + p \frac{c_0}{\l_j^{(n-2)/2}}  \int_{\R^n_+} \d_i ^{p-1}{\l _i} \frac{\partial\d_i }{\partial \l _i } H(a_j,.) + O\Big( \frac{1}{(\l_jd_j)^2} \int \d_i^p\d_j\Big).\ee

The first term is estimated in \eqref{qq11}. For the second one, expanding $H(a_j,.)$ around $a_i$ in the ball $B(a_i,d_i/2)$, it follows
\begin{align*} \frac{p c_0}{\l_j ^{(n-2)/2}}  \int_{\R^n_+} \d_i ^{p-1} & {\l _i} \frac{\partial\d_i }{\partial \l _i } H(a_j,.)  =  \frac{H(a_j,a_i)}{(\l_j\l_i) ^{(n-2)/2}} pc_0^{\frac{2n}{n-2}}\frac{n-2}{2} \int_{B(0,\l_i d_i/2)}\frac{ 1-|x|^2}{(1+|x|^2)^{(n+4)/2}} \\
&  + O\Big( \frac{1}{ d_i^2} \int _{B(a_i,d_i/2)} |x-a_i|^2 \d_i ^p \d_j + \int _{\R^n\setminus B(a_i,d_i/2)} \d_i ^p \d_j \Big)\\
& = -\frac{n-2}{2} c_2 \frac{H(a_j,a_i)}{(\l_j\l_i) ^{(n-2)/2}} + O\Big(  \frac{1}{(\l_i d_i)^2} \e_{ij} (\ln \e _{ij }^{-1}  )^{\frac{n-2}{n}} + \frac{1}{(\l_j d_j)^n} + \frac{1}{(\l_i d_i)^n} \Big).
\end{align*}
For $j = i$, using again Lemma  \ref{lem:varphi} and expanding $H(a_i,.)$ around $a_i$ in the ball $B(a_i,d_i/2)$, it follows
\begin{align*}  \langle \varphi_i,{\l _i} \frac{\partial\varphi _i }{\partial \l _i } \rangle & = p  \int_{\R^n_+} \d_i ^{p-1}{\l _i} \frac{\partial\d_i }{\partial \l _i } \d_i + p \frac{c_0}{\l_i^{(n-2)/2}}  \int_{\R^n_+} \d_i ^{p-1}{\l _i} \frac{\partial\d_i }{\partial \l _i } H(a_i,.) + O\Big( \frac{1}{\l_i ^{(n+2)/2}d_i^n} \int \d_i^p\Big)\\
& = - p\int_ {\R^n_- }  \d_i ^{p}{\l _i} \frac{\partial\d_i }{\partial \l _i } + p \frac{c_0^{2n/(n-2)}}{\l_i^{n-2}}\frac{n-2}{2} H(a_i,a_i) \int_{B(0,\l_i d_i/2)}\frac{ 1-|x|^2}{(1+|x|^2)^{(n+4)/2}} \\
& + O\Big( \frac{1}{\l_i^{(n-2)/2} d_i ^n} \int _{B(a_i,d_i/2)} | x-a_i |^2 \d_i^p  +  \frac{1}{\l_i^{(n-2)/2} d_i ^{n-2}}\int_{\R^n \setminus B(a_i,d_i/2)} \d_i^p  +  \frac{1}{(\l_i d_i )^{n}}\Big) \\
& =  - \frac{n-2}{2} c_2 \frac{H(a_i,a_i)}{\l_i ^{n-2}}  + O\Big( \frac{\ln(\l_i d_i)}{(\l_i d_i)^n} \Big).
\end{align*}
In addition, the analogue of \eqref{aaz4} becomes
\begin{align*}
 \langle \o , \l_i \frac{\partial \varphi_i }{ \partial \l_i}  \rangle & = (\i \o)({a}_i)\,  p \int_{\R^n } \d_i^{p-1} \l_i \frac{\partial \d_i }{ \partial \l_i} + O\Big( \int_{\R^n \setminus B(a_i, d_i)} \d_i ^p + \int_{ B(a_i, d_i) } \d_i ^p | x- {a}_i|^2 \Big)\\
 & = - 2\ov{c}_6 \frac{\o (a_i) }{\l_i ^{(n-2)/2}} + O\Big( \frac{\ln(\l_i d_i)}{\l_i ^{(n+2)/2}} + \frac{1}{(\l_i d_i )^n}\Big). \end{align*}
Hence the first part of \eqref{nablaI} follows. For the second part of \eqref{nablaI}, observe that \eqref{az5}, \eqref{az6} and \eqref{wdw2} hold with $\varphi_i$ instead of $\d_i$.
We remark that the first term of the right hand side of \eqref{az6} has to be changed and we have
$$  \int_{\mathbb{S}^n_+ } K \varphi_i  ^{ p-1-\e} \ov{v}  \l_i \frac{ \partial \varphi_i }{ \partial \l_i} =  \int_{\mathbb{S}^n_+ } K \d_i  ^{ p-1-\e} \ov{v}  \l_i \frac{ \partial \d_i }{ \partial \l_i} + O\Big( \frac{1}{ (\l_i d_i ^2)^{(n-2)/2}} \int_{ B(a_i, d_i) } \d_i ^{p-1} | \ov{v} | +  \int_{ \mathbb{S}^n_+\setminus B(a_i, d_i) } \d_i ^{p} | \ov{v} | \Big)$$ and the estimate follows from \eqref{az7}, \eqref{2***} and \eqref{1***} (by using the Holder's inequality).

Therefore \eqref{fgh1} holds true with $\varphi_i$ instead of $\d_i$. It remains to estimate the analogue of  \eqref{wdw2}. Observe that, for $j \leq q$, it holds
\begin{align} \int_{\mathbb{S}^n_+ } K  \d_j ^{p-\e} \l_i \frac{\partial \varphi_i }{ \partial \l_i} = & \frac{ K(a_j) c_0^{-\e}}{\l_j^{\e(n-2)/2}} \langle  \d_j , \l_i \frac{\partial \varphi_i }{ \partial \l_i} \rangle + O\Big( \int_{\mathbb{S}^n_+  } d( x,a_j) \d_j^{p-1} (\d_j\d_i) \Big)\\
& +O \Big( \e  \int_{\mathbb{S}^n_+ } \d_j ^{p-1}\ln [ 2+(\l_j^2-1)(1-\cos d (x,a_j)) ] (\d_j \d_i)  \Big).\nonumber
\end{align}
These terms are computed in  \eqref{qq11} and \eqref{wxcv1}. Now, for $j \geq q+1$ with $j\neq i$, it holds
\begin{align} \int_{\mathbb{S}^n_+ } K  \varphi_j ^{p-\e} \l_i \frac{\partial \varphi_i }{ \partial \l_i} = &  \int_{\mathbb{S}^n_+ } K  \d_j ^{p-\e} \l_i \frac{\partial \varphi_i }{ \partial \l_i}  + O\Big( \int_{\mathbb{S}^n_+ } \d_j^{p-1} | \varphi_j - \d_j | \d_i \Big)\label{qq1122}\\
& = \frac{ K(a_j) c_0^{-\e}}{\l_j^{\e(n-2)/2}} \langle  \varphi_j , \l_i \frac{\partial \varphi_i }{ \partial \l_i} \rangle + O\Big(  \int_{\mathbb{S}^n_+ } \d_j^{p-1} | \varphi_j - \d_j | \d_i  +  \int_{\mathbb{S}^n_+ } d(x,a_j) \d_j^{p-1}(\d_j  \d_i) \nonumber\\
& + \e  \int_{\mathbb{S}^n_+ } \d_j ^{p-1}\ln [ 2+(\l_j^2-1)(1-\cos d (x,a_j)) ] (\d_j \d_i)  \Big) \nonumber
\end{align}
and these terms are computed in  \eqref{qq112} and \eqref{wxcv1}.  Finally, it remains the case $j = i$.
Note that, since $a_i \in  \mathbb{S}^n_+$, it follows that $B_i:= B(a_i,d_i/2) \subset  \mathbb{S}^n_+$ and $\l_i d_i$ is very large. Using Lemma \ref{lem:varphi}, it follows that
\begin{align}
& \int_{\mathbb{S}^n_+ } K  \varphi_i ^{p-\e} \l_i \frac{\partial \varphi_i }{ \partial \l_i} \label{jsec} \\
 & = \int_{B_i} K \d_i^{p-\e}  \l _i \frac{\partial\d_i }{\partial\l _i} - \frac{n-2}{2} \int_{B_i} K \d_i^{p-\e} \frac{ c_0 H(a_i,.)}{\l_i^{(n-2)/2}} + (p-\e) \int_{B_i} K \d_i^{p-\e-1}  \l _i \frac{\partial\d_i }{\partial\l _i}  \frac{ c_0 H(a_i,.)}{\l_i^{(n-2)/2}} \nonumber \\
& + O \Big( \int_{B_i} \d_i^p \frac{1}{\l_i^{(n+2)/2} d_i^n} +  \int_{B_i} \d_i^{p-1} \frac{ H(a_i,.)^2}{\l_i^{n-2}} + \int_{B_i} \d_i^{p-1} \frac{1}{\l_i^{(n+2)/2} d_i^n} \frac{ H(a_i,.)}{\l_i^{(n-2)/2}} + \int_{ \mathbb{S}^n_+\setminus B_i} \d_i ^{p+1} \Big).\nonumber \end{align}
The estimate of the last term is given in \eqref{2***} and, using Lemma \ref{lem:varphi}, easy computations imply that
\begin{align*}
& \int_{B_i} \d_i^{p-1} \frac{1}{\l_i^{(n+2)/2} d_i^n} \frac{ H(a_i,.)}{\l_i^{(n-2)/2}} \leq \frac{c}{\l_i^{(n+2)/2} d_i^n}  \int_{B_i} \d_i^{p} \leq \frac{c}{(\l_i d_i)^n} , \\
&   \int_{B_i} \d_i^{p-1} \frac{ H(a_i,.)^2}{\l_i^{n-2}} \leq \frac{c}{ \l_i^{ n-2} d_i^{2n-4}} \int _{B_i} \d_i^{p-1} \leq \frac{c}{(\l_i d_i)^n} .
\end{align*}
Furthermore, using Lemma \ref{lowerL2}, the second and the third terms of \eqref{jsec} follow from the following estimates
$$  \e \int_{B_i} K \d_i^{p} \ln [2+(\l_i ^2-1)(1-\cos d(a_i,x) ) ]  \frac{ c_0 H(a_i,.)}{\l_i^{(n-2)/2}}   \leq  \frac{ \e}{( \l_i d_i)^{n-2} }, $$
\begin{align*} \int_{B_i} K \d_i^{p} \frac{ c_0 H(a_i,.)}{\l_i^{(n-2)/2}} & = \frac{ c_0 H(a_i,a_i) K(a_i)}{\l_i^{(n-2)/2}} \int_{B_i}  \d_i^{p}  + O \Big( \frac{ c}{ \l_i^{(n-2)/2} d_i ^n } \int_{B_i} | x- a_i |^2 \d_i^{p} \Big) \\
& = c_2K(a_i) \frac{  H(a_i,a_i) }{\l_i^{n-2}} + O\Big( \frac{\ln(\l_i d_i)}{ (\l_i d_i )^n} \Big) \qquad (\mbox{using \eqref{3***} and \eqref{2***} }) \end{align*}
\begin{align*}
p \int_{B_i} K \d_i^{p-1}  \l _i \frac{\partial\d_i }{\partial\l _i}  \frac{ c_0 H(a_i,.)}{\l_i^{(n-2)/2}} & =   \frac{ c_0 H(a_i,a_i) K(a_i)}{\l_i^{(n-2)/2}} p \int_{B_i}  \d_i^{p-1}  \l _i \frac{\partial\d_i }{\partial\l _i} +O \Big( \frac{ c}{ \l_i^{(n-2)/2} d_i ^n } \int_{B_i} | x- a_i |^2 \d_i^{p} \Big) \\
& = - c_2 \frac{n-2}{2} K(a_i) \frac{  H(a_i,a_i) }{\l_i^{n-2}} + O\Big( \frac{\ln(\l_i d_i)}{ (\l_i d_i )^n} \Big) \quad (\mbox{using \eqref{4***} and \eqref{2***}})
\end{align*}
where we have used for the first one the fact that $H(a_i, .) \leq c /  d_i ^{n-2} $ and in the other, we have expanded $K H(a_i,.)$ around $a_i$.\\
It remains to estimate the first integral of \eqref{jsec}. Expanding as in \eqref{dd1}, we need to compute
\begin{align}
 \int_{B_i} {K} \d_i^{p}  \l _i \frac{\partial\d_i }{\partial\l _i} & =  \int_{\tilde{B}_i} \wtilde{K} \d_i^{p}  \l _i \frac{\partial\d_i }{\partial\l _i}  \nonumber  \\
 & = O\Big( \int _{\R^n\setminus \tilde{B}_i} \d_i ^{p+1}\Big) + \frac{1}{2} \sum_k \frac{\partial ^2\wtilde{K}}{\partial x_k^2} (a_i) \int _{\tilde{B}_i} (x-a_i)_k^2 \d_i^{p}  \l _i \frac{\partial\d_i }{\partial\l _i} + O\Big( \int _{\tilde{B}_i} | x-a_i|^4 \d_i^{p+1} \Big) \nonumber \\
  & = - \frac{c_7}{4} \frac{\D \wtilde{K} (a_i)}{\l_i ^2} + O\Big( \frac{1}{\l_i ^4} + \frac{1}{(\l_i d_i )^n}\Big) \qquad \mbox{ where } c_7:= \frac{n-2}{n} c_0^{p+1} \int_{\R^n} \frac{ | x |^2( |x|^2 - 1 )}{(1+| x |^2)^{n+1}}\label{c7}
   \end{align}

\begin{align*}  \int_{B_i} { K} \d_i ^p & \ln [ 2+(\l_i ^2-1)  (1-\cos(d(a_i,x))) ]  \l _i \frac{\partial\d _i }{\partial\l _i}  =  \int_{\tilde{B}_i} \wtilde{ K} \d_i ^p \ln\Big(\frac{2}{1+|y|^2}(1+\l_i ^2 | y| ^2)\Big) \l _i \frac{\partial\d _i }{\partial\l _i} \\
&= K(a_i) c_0^{p+1} \frac{n-2}{2} \int_{\R^n}   \frac{ 1 - | x |^2}{(1+| x |^2)^{n+1}} \ln(1+|x|^2)+ O\Big(  \frac{1}{\l_i ^2} +\frac{\ln(\l_i d_i)}{ (\l_i d_i )^n} \Big) \\
& = - 2  K(a_i)  \frac{2 c_4}{n-2} +  O\Big(   \frac{1}{\l_i ^2} +  \frac{\ln(\l_i d_i)}{ (\l_i d_i )^n} \Big) \end{align*}
where $\tilde{B}_i$ is the image of $B_i$ by the change of variables.\\
Combining the previous estimates the proof of the third claim follows.

The last claim follows as the previous one by noting that \eqref{qq11}-\eqref{qq1122} are seen as $O(\e_{ij})$. We omit the details.
\end{pf}

\subsection{Expansion in the neighborhood at infinity $V(q,\ell,\tau)$ }

In this subsection we write down, for later use,  the expansion of the gradient in case $\o = 0$. These expansions follow from the previous expansions  by simply  taking $\o =0$
\begin{pro}\label{devIeps}
Let   $u:= \sum_{i=1}^{N}\a_i\varphi_{a_i,\l_i}  + \ov{v}\in V(q,\ell,\tau) $ (with $N:= q+\ell$). Assume that $\e \ln \l_i$ is small for each $i$. For each $i \leq q$, there hold
\begin{align*}
& \langle \n  I_\e(u),\d_i\rangle = \a_iS_n \, (1-\a_i^{p-1}K(a_i) \l_i^{-\e(n-2)/2} )+ O\Big(\e  + \sum \e_{ij} +   \frac{ 1}{\l_i } +  R(\e,a,\l)^2\Big),\\
 & \langle \n I_\e(u),\l_i\frac{\partial\d_i}{\partial \l_i}\rangle = \frac{\a_i^{p}}{\l_i^{\e(n-2)/2}} \Big( c_4 K(a_i) \, \e - \frac{c_3}{\l_i}\frac{\partial K}{\partial \nu}(a_i) \Big) + O\big(\e^2  + \sum_{j \neq i}\e_{ij} + \frac{1}{\l_i^2}  +  R^2(\e,a,\l) \big),\\
 &  \langle \n I_\e(u),\frac{1}{\l_i}\frac{\partial\d_i}{\partial a_i}\rangle _{\lfloor T_{a_i}(\partial \mathbb{S}^n_+)}  =  \frac{c_2}{2}  \sum_{i \neq j \leq q}\a_j \frac{1}{\l_i}\frac{\partial \e_{ij}}{\partial a_i} \Big(1- \frac{\a_i^{p-1}K(a_i)}{ \l_i^{\e(n-2)/2}}  - \frac{\a_j^{p-1}K(a_j)}{ \l_j^{\e(n-2)/2}} \Big)\\
& \qquad \qquad\qquad  \qquad\qquad - \frac{\a_i^{p} }{ \l_i^{\e(n-2)/2}} c_5\frac{\n_T K(a_i)}{\l_i} + O\Big( \frac{1}{\l_i ^2  } + R_{3,i} (\e,a,\l)    \Big)
\end{align*}
where $T_{a_1}(\partial \mathbb{S}^n_+)$ denotes  the tangent space at  $a_i$ and the parameters are defined in Proposition \ref{devIeps3}.
\end{pro}

\begin{pro}\label{devIeps2}
Let   $u:= \sum_{i=1}^{N}\a_i\varphi_{a_i,\l_i}  + \ov{v}\in V(q,\ell,\tau) $ (with $N:= q+\ell$). Assume that $d_i \geq c > 0$ for each $i > q$ and $\e \ln \l_i$ is small for each $i$. Then, for each $i > q$, there hold

\begin{align*}
& \langle \n  I_\e(u),\varphi_i\rangle = 2 \a_iS_n \, \Big(1- \frac{ \a_i^{p-1}K(a_i) }{ \l_i^{\e(n-2)/2} } \Big) + O\Big(\e  + \sum \e_{ij} +   \frac{ 1}{\l_i ^2} +  R^2(\e,a,\l) + \frac{1}{(\l_i d_i)^{n-2}}  \Big),\\
 & \langle \n I_\e(u),\l_i\frac{\partial\varphi_i}{\partial \l_i}\rangle =   \frac{\a_i^{p}}{\l_i^{\e(n-2)/2}} \Big(2 c_4 K(a_i) \, \e + {c_7}\frac{\D K(a_i)}{\l_i^2} \Big) - {c_2}\frac{n-2}{2} {\a_i} \frac{ H(a_i,a_i)}{ \l_i ^{n-2}} \times \\
 & \qquad \qquad  \times   \Big( 1 - 2 \frac{\a_i^{p-1}K(a_i)}{ \l_i ^{\e(n-2)/2}} \Big)+ O\Big(\e^2  + \sum_{j \neq i}\e_{ij}   +  R^2(\e,a,\l)  + \sum_{ j > q} \frac{\ln(\l_k d_k)}{(\l_k d_k)^n}\Big),\\
 &  \langle \n I_\e(u),\frac{1}{\l_i}\frac{\partial \varphi_i}{\partial a_i}\rangle =  - 2\a_i^{p} \l_i^{-\e(n-2)/2} c_5\frac{\n K(a_i)}{\l_i} +  O\Big(\e^2  + \sum_{j \neq i}\e_{ij} + \frac{1}{\l_i^2}  + R^2(\e,a,\l)\Big).
\end{align*}
\end{pro}

\section{Critical points of the Kirchoff-Routh type functional $\mathcal{F}_{z,m}$}

In this section we study the existence of critical points of the Kirchoff-Routh type functional $\mathcal{F}_{z,m}$. Our main result states as follows

\begin{pro}\label{pF}
Let $z \in \partial \mathbb{S}^n_+$ be a non degenerate critical point of $K_1:= K_{\lfloor{\partial \mathbb{S}^n_+}}$ and $\mathcal{F}_{z,m}$ be defined in \eqref{dF2}. Then we have
\begin{enumerate}
\item[(i)] If $z$ is a local minimum then $\mathcal{F}_{z,m}$ achieves its  minimum.
\item[(ii)]  If $z$ is a local maximum  then $\mathcal{F}_{z,m}$ does not have any critical point.
\item[iii)]
If the hessian matrix $D^2K_1(z)$ has at least a positive eigenvalue $\l$, then  the function $\mathcal{F}_{z,2}$ admits at least one critical point.
\item[iv)]
For each simple positive eigenvalue $\l$ of   the hessian matrix $D^2K_1(z)$  the function $\mathcal{F}_{z,2}$  has a  non-degenerate critical point $(x_\l,y_\l)$.\\
Furthermore, if $\l_1$ and $\l_2$ are two simple positive eigenvalues, then $(x_{\l_1},y_{\l_1}) \neq (x_{\l_2},y_{\l_2})$.
\end{enumerate}
\end{pro}

\begin{pf} The proof of $(i)$ is trivial  since  $D^2K_1(z)$ is a positive definite matrix. Indeed  in this case $\mathcal{F}_{z,m}$ is coercive.\\
The proof of statement $(ii)$ is as follows:  writing $(\xi_1,\cdots,\xi_m):= r \Lambda$ with $r >0$ and $\L := (\L_1,\cdots,\L_m) \in (\R^{n-1})^m$ with $\| \L \| =1$, the function $\mathcal{F}_{z,m}$  reads  as follows:
$$\mathcal{F}_{z,m}(r, \L):= \frac{1}{2} r ^2 \sum_{i=1}^{m} D^2K_1(z)( \L_i , \L_i ) + \frac{1}{r^{n-2}} \sum_{1 \leq i < j \leq m}\frac{1}{| \L_i - \L_j | ^{n-2}} .$$
Hence, under the assumption that  $D^2K_1(z)( \L_i , \L_i ) \leq 0$ for each $\L_i$, we have
$$\frac{\partial \mathcal{F}_{z,m}}{\partial r} = r  \sum_{i=1}^{m} D^2K_1(z)( \L_i , \L_i ) - \frac{n-2}{r^{n-1}} \sum_{1 \leq i < j \leq m}\frac{1}{| \L_i - \L_j | ^{n-2}} < 0 $$
Hence the function $\mathcal{F}_{z,m}$ does not have any critical point if the point $z$ is a local maximum. \\
To prove the statement $(iii)$ let $\l$ be a positive eigenvalue of $M:= D^2K_1(z)$ (which exists by our assumption). To find  a critical point of $\mathcal{F}_{z,2}$, we need to solve the equation $D \mathcal{F}_{z,2} (x,y)=0$ which is equivalent to
\be\label{e:F1}
 M.x -  (n-2)  \frac{x-y}{| x - y | ^{n}} = 0 \quad \mbox{ and } \quad
 M.y -  (n-2)  \frac{y-x}{| x - y | ^{n}} = 0.
\ee
Summing the two equations, we get $M.(x+y)=0$ which implies that $x=-y$ (since $M$ is assumed to be non-degenerate). Inserting  this information in the first equation of \eqref{e:F1}, we obtain
\be\label{e:F2} M.x =   \frac{(n-2)}{2^{n-1}}  \frac{x}{| x  | ^{n}} := {\underline{c}}   \frac{x}{| x  | ^{n}} . \ee
It is easy to see that  $ {\ov{x}} := ( {\underline{c}} /\l)^{1/n} u_\l $ is a solution of \eqref{e:F2}  (where $u_\l$ is an  unit eigenvector corresponding  to the eigenvalue $\l$). Hence the proof of the first claim.\\
For the last claim, we have by assumption  that the eigenspace  $ V_\l $ of  the eigenvalue $\l$  is one dimensional and our goal is  to prove that the critical point $({\ov{x}} , -{\ov{x}} )$ obtained  in the statement $(iii)$ of this proposition is non-degenerate.  To that aim, we observe that

\begin{equation*}
D^2 \mathcal{F}_{z,2}({\ov{x}} ,{\ov{y}} ) ((x,y),.) =
  \begin{pmatrix}
    M.x - {(n-2)}  \frac{x-y}{| \ov{x} - \ov{y}  | ^{n}} + n (n-2)    \frac{\langle \ov{x}-\ov{y}, x-y\rangle }{| \ov{x} - \ov{y}  | ^{n+2}} ( \ov{x}-\ov{y}) \\
    M.y -   (n-2)  \frac{ y - x}{|  \ov{y} - \ov{x}  | ^{n}} + n (n-2)  \frac{\langle \ov{y}-\ov{x}, y-x\rangle }{| \ov{x} - \ov{y}  | ^{n+2}} ( \ov{y}-\ov{x})
  \end{pmatrix} .
\end{equation*}

We claim that  the unique solution of  $D^2 \mathcal{F}_{z,2}({\ov{x}} ,{\ov{y}} ) ((x,y),.) = 0$  is $(0,0)$. Indeed, summing up the two components, we obtain  that $M.(x+y)=0$ which implies that $x=-y$. The first component hence becomes:
\be\label{e:F3} M.x =  \frac{(n-2)}{2^{n-1}}   \frac{x}{| \ov{x}  | ^{n}} -  \frac{n(n-2)}{2^{n-1}}   \frac{\langle \ov{x}, x \rangle }{| \ov{x}  | ^{n+2}} \ov{x}. \ee
Taking the scalar product of \eqref{e:F3} with $\ov{x}$, we get
$$ \langle M. x , \ov{x} \rangle =  \frac{(n-2)}{2^{n-1}}  \frac{ \langle x, \ov{x} \rangle} {| \ov{x}  | ^{n}} - \frac{n(n-2)}{2^{n-1}}   \frac{\langle \ov{x}, x \rangle }{| \ov{x}  | ^{n}} =   \frac{(n-2)}{2^{n-1}}(1 - n) \frac{ \langle x, \ov{x} \rangle} {| \ov{x}  | ^{n}} . $$
Using the fact that $M$ is a symmetric matrix, we derive that $  \langle M. x , \ov{x} \rangle = \langle M. \ov{x} , {x} \rangle $. Using the fact that $\ov{x}$ is a solution of \eqref{e:F2}, the previous equation implies
$$   \frac{(n-2)}{2^{n-1}}(1- n) \frac{ \langle x, \ov{x} \rangle} {| \ov{x}  | ^{n}} =  \frac{(n-2)}{2^{n-1}}  \frac{\langle \ov{x}, x \rangle }{| \ov{x}  | ^{n}}$$ which implies that $\langle x, \ov{x} \rangle = 0$. Putting this information in \eqref{e:F3} and using the fact that $| \ov{x} | ^n = \underline{c} /\l$ we derive that $M. x = \l x$ which implies that $x \in V_\l$. Now, since we assumed that $\mbox{dim }V_\l=1$ and using the fact that $\langle x, \ov{x} \rangle = 0$ and $\ov{x} \in V_\l$ we obtain that $x=0$ and therefore $(x,y)=(0,0)$. Hence $D^2 F({\ov{x}} ,{\ov{y}} ) $ is non-degenerate. This completes the proof.
\end{pf}

\begin{pro} \label{p:dF234}
Let $m\in \N$ and let $z \in \partial \mathbb{S}^n_+$ be a non degenerate critical point of $K_1:= K_{\lfloor{\partial \mathbb{S}^n_+}}$ and $\mathcal{F}_{z,m}$ be defined in \eqref{dF2}. Let $M$ be the associated matrix to $D^2 K_1(z)$. We assume that $M$ has one positive simple eigenvalue and the others are negative. Then  $\mathcal{F}_{z,m}$ admits a non degenerate critical point.
\end{pro}
\begin{pf}
Let $\s_k$ (for $k=1,\cdots,n-1$) be the eigenvalues of $M$ (not necessary different). By our assumption, we have $ \s_1 > 0$ and $\s_k < 0$ for each $k \geq 2$. Let $(e_1,\cdots e_{n-1})$ be an orthonormal basis associated to the eigenvalues $\s_k$. In a first step, we define 
 $f_m$ to be 
$$f_m(\xi) = \frac{1}{2} \s_1 | \xi |^2 + \sum _{1\leq i < j \leq m} \frac{1}{| \xi_i - \xi_j |^{n-2}} \quad  \mbox{ where } \xi := (\xi_1,\cdots,\xi_m) \in \R^{m}\setminus \G$$
with $\G := \{ \xi \in \R^m: \exists \, i\neq j \mbox{ with } \xi_i = \xi_j\}$.\\
Observe that $f_m$ achieves its infimum (since $\s_1 > 0$). Let $\ov{\xi}$ be such that  $f_m(\ov{\xi}) = \inf f_m$. Hence $\ov{\xi}$ is a critical point of $f_m$. Note that, for $X =(x_1,\cdots,x_m) \in \R^m$, easy computations imply that
$$D^2 f_m(\ov{\xi}) (X,X) = \s_1 | X |^2 + (n-2)(n-1)  \sum_{1\leq i < j \leq m} \frac{ | x_i - x_j |^2 }{| \ov{\xi}_i - \ov{\xi} _j |^{n} }.$$
Since $\s_1 > 0$, it follows that $\ov{\xi}$ is a non degenerate critical point of $f_m$.\\
Now, let $ \ov{Y} := (\ov{Y}_1, \cdots,\ov{Y}_m)$ with $\ov{Y}_k= \ov{\xi}_k e_1$ for each $k \geq 1$. We claim that\\
{\bf CLAIM: }  $ \ov{Y} $ is a non degenerate critical point of $\mathcal{F}_{z,m}$.\\
In fact, let $X:= (X_1, \cdots,X_m)$ with $X_k:= (X_{k,1},\cdots,X_{k,{n-1}})\in \R^{n-1}$ for each $k$. It holds that

\begin{align*} D \mathcal{F}_{z,m} (\ov{Y}) (X) & = \sum_{i=1}^m \langle M \ov{Y}_i, X_i \rangle - (n-2) \sum _{1\leq i < j \leq m} \frac{ \langle \ov{Y}_i -  \ov{Y}_j, X_i - X_j \rangle }{ |  \ov{Y}_i -  \ov{Y}_j |^n} \\
& = \sum_{i=1}^m \ov{\xi}_i \s_1 X_{i,1} - (n-2) \sum _{1\leq i < j \leq m} \frac{ ( \ov{\xi}_i -  \ov{\xi}_j) (X_{i,1} - X_{j,1} ) }{ |  \ov{\xi}_i -  \ov{\xi}_j |^n} \\
& = D f_m(\ov{\xi}) (X^1) = 0  \quad \forall \, \, X \in (\R^{n-1})^m \quad \mbox{ where } X^1 := (X_{1,1},\cdots,X_{m,1}).
\end{align*}
Hence $\ov{Y}$ is a critical point of $\mathcal{F}_{z,m}$. Now we need to solve the following equation: $ D^2 \mathcal{F}_{z,m} (\ov{Y}) (X,.) = 0$ which is equivalent to
\begin{equation*}
M X_i = (n-2) \sum_{j \neq i} \frac{ X_i - X_j }{| \ov{Y}_i - \ov{Y}_j | ^n} - n(n-2) \sum_{j \neq i} \frac{ \langle \ov{Y}_i -  \ov{Y}_j, X_i - X_j \rangle  }{| \ov{Y}_i - \ov{Y}_j | ^{n+2}} (\ov{Y}_i - \ov{Y}_j ) \quad \mbox{ for each } i = 1,\cdots,m.
\end{equation*}
Recall that $ \ov{Y}_i = \ov{\xi}_i e_1$ and let $X_i:= \sum_{k=1}^{n-1} X_{i,k} e_k$. Hence the previous equation becomes
\begin{equation*}
\sum_{k=1}^{n-1} X_{i,k} \s_k e_k = (n-2) \sum_{j \neq i} \frac{ 1 }{| \ov{\xi}_i - \ov{\xi}_j | ^n} \sum_{k=1}^{n-1} (X_{i,k} - X_{j,k} )e_k  - n(n-2) \sum_{j \neq i} \frac{ X_{i,1} - X_{j,1}  }{| \ov{\xi}_i - \ov{\xi}_j | ^{n}} e_1
\end{equation*}
which implies that
\begin{align}
& \s_1 X_{i,1} = (n-2)(1-n) \sum_{j \neq i} \frac{ X_{i,1} - X_{j,1}  }{| \ov{\xi}_i - \ov{\xi}_j | ^n}  \quad  \quad \mbox{ for each } i = 1,\cdots,m  \mbox{ and } \label{0909} \\
& \s_k X_{i,k} = (n-2) \sum_{j \neq i} \frac{ X_{i,k} - X_{j,k}  }{| \ov{\xi}_i - \ov{\xi}_j | ^n}  \quad \mbox{ for each } k \geq 2  \mbox{ and for each } i = 1,\cdots,m.  \label{0908}
\end{align}
 Multiplying \eqref{0909} by $X_{i,1}$ and summing over $i = 1,\cdots,m$, we obtain
 $$ \s_1\sum_{i=1}^m X_{i,1}^2 = (n-2)(1-n)    \sum_{1 \leq i < j \leq m} \frac{ (X_{i,1} - X_{j,1} )^2 }{| \ov{\xi}_i - \ov{\xi}_j | ^n}  \leq 0 $$ which implies that $X_{i,1} = 0$ for ech $i=1,\cdots,m$ (since $\s_1 >0$).\\
In the same way, for $k \geq 2$, multiplying \eqref{0908} by $X_{i,k}$ and summing over $i = 1,\cdots,m$, we obtain
 $$ \s_k\sum_{i=1}^m X_{i,k}^2 = (n-2)   \sum_{1 \leq i < j \leq m} \frac{ (X_{i,k} - X_{j,k} )^2 }{| \ov{\xi}_i - \ov{\xi}_j | ^n}  \geq 0 $$ which implies that $X_{i,k} = 0$ for ech $i=1,\cdots,m$ (since $\s_k  <0$).\\
 Hence it follows that the unique solution of $ D^2 \mathcal{F}_{z} (\ov{Y}) (X,.) = 0$ is $ X = 0$ which implies that the critical point $\ov{Y}$ is non degenerate. Therefore the proof of the proposition is completed
\end{pf}

\section{Proof of the main results}

The strategy of the proofs of the theorems is the following: In each proof, we will define a set $M_\e$ (depending on the kind of the blow up points which we need to obtain). Concerning Theorems \ref{t:13}, \ref{t:11} and \ref{t:12} (the zero limit case), the elements of $M_\e$ are some points $\mathcal{M} :=(\a,\l,x,v) \in (\R_+)^N\times  (\R_+)^N\times (\ov{ \mathbb{S}}^n_+)^N \times H^1(\mathbb{S}^n_+)$  where $v$ satisfies \eqref{V0} and the other variables satisfy some conditions and where $N$ is the number of the bubbles in the desired constructed solution. In addition, we will define a function
\begin{eqnarray}\label{F10}  \Psi_{\e}:M_{\e}\,\rightarrow \R ;\quad \mathcal{M}=(\a,\l,x,v)\mapsto I_\e\big( \sum_{i=1}^{N} \a_i \varphi_{(x_i,\l_i)} \, + \, v\big). \end{eqnarray}
The first step is inspired from \cite{BLR}. Since $v$ satisfies some orthogonality conditions, using the Euler-Lagrange's coefficients, it is easy to get the following proposition.

\begin{pro}\label{G1} Let $ \mathcal{M} =(\a,\l,x,v)\in M_{\e}$. $\mathcal{M}$ is a critical point of $\Psi_{\e}$ if and only if $ \sum_{i=1}^{N}\a_i\varphi_{(x_i,\l_i)}   +  v$ is a critical point of $I_\e$, i.e. if and only if there exists $\big(A,B,C\big)\in \R^N\times\R^N\times\big(\R^{n}\big)^N$ such that the following holds :
\begin{align}
 & (E_{\a_i}) \qquad \frac{\partial \Psi_{\e}}{\partial \a_i}(\a,\l,x,v)=0,\,\, \forall \,i=1,\cdots,N \notag\\
 & (E_{\l_i}) \qquad \frac{\partial \Psi_{\e}}{\partial \l_i}(\a,\l,x,v)=B_i\langle \l_i \frac{\partial^2 \varphi_i}{\partial\l_i^2},v\rangle+\sum_{j=1}^{n} C_{ij}\langle \frac{1}{\l_i}\frac{\partial^2 \varphi_i}{\partial x_i^j\partial\l_i},v\rangle,\,\ \forall \,i=1,\cdots, N\notag\\
 & (E_{x_i}) \qquad \frac{\partial \Psi_{\e}}{\partial x_i}(\a,\l,x,v)\lfloor_{T_{x_i}( \mathbb{S}^n_+)} \, = \, B_i\langle \l_i \frac{\partial^2 \varphi_i}{\partial\l_i\partial x_i},v\rangle+ \sum_{j=1}^{n} C_{ij} \langle \frac{1}{\l_i}\frac{\partial^2 \varphi_i}{\partial x_i^j\partial x_i},v\rangle, \,\ \forall \,i=1,\cdots,N \notag\\
 & (E_v)\qquad \frac{\partial \Psi_{\e}}{\partial v}(\a,\l,x,v)=\sum_{i=1}^N \biggl(A_i \varphi_i+B_i\l_i \frac{\partial \varphi_i}{\partial \l_i}+\sum_{j=1}^{n} C_{ij}\frac{1}{\l_i}\frac{\partial \varphi_i}{\partial x_i^j}\biggr)\notag
\end{align}
where $\varphi_i:=\varphi_{(x_i,\l_i)} $ and  where $C_{in}=0$ if $a_i \in \partial \mathbb{S}^n_+$.
\end{pro}
The results of the theorems will be obtained through a careful analysis of the previous equations on $M_{\e}$. In the following we will present the proof of each theorem and give the precise  definition of the sets involved  in each case.

 \subsection*{ Proof of Theorem \ref{t:11} }

 For $z$ a non-degenerate critical point of $K_1$ (not a local maximum) and $2 \leq m \in \N$, we define the following set:
\begin{align}\label{Meps}  M_\e(z) :=M_{\e}(z,\mu,m):=  \{ & \mathcal{M} :=  (\a,\l,x,v)\in(\R_+)^m \times(\R_+)^m\times (\partial \mathbb{S}^n_+)^m \times H^1(\mathbb{S}^n_+):  \\
& v\mbox{ satisfies \eqref{V0}} ; \, \,   |{\a_i^{4/(n-2)}K(x_i)} - 1 | < \e \ln^2\e\, ; \, \, C^{-1} \e \leq \l_i^{-1} \leq  C \e ; \nonumber\\
& d(x_i , z )  \leq \mu \,\, \forall \, i\, ;  \, {C}^{-1} \e^{(n-2)/n} \leq d( x_i , x_ j ) \leq C\, \e^{(n-2)/n}\quad \forall i \neq j  \}, \nonumber
\end{align}
where $C$ is a large positive constant and $\mu$ is a small positive constant.\\
We notice that, for $(\a,\l,x,0)\in M_\e(z,\mu,m)$, one has,  for $i \neq j$, that $d( x_i , x_j ) \to 0$ as $\e \to 0$ and standard computation implies that
$ 2 \big( 1 - \cos d(x_i , x_j ) \big) \sim d( x_i , x_j )^2 $  and therefore, using \eqref{epsilon}, we get
\be\label{eps1} \e_{ij} = O (\e ^{1+ (n-4)/n}).\ee
Next, as said before, on $M_\e$,   we  define the function $  \Psi_{\e} $ (introduced in \eqref{F10}). Note that Proposition \ref{G1} holds with $N=m$ and in our case, $\varphi_{a_i,\l_i} = \d_{a_i,\l_i} $ for each $i$.\\
In the sequel and for   the  sake of simplicity, we will write $\d_i$ and $\varphi_i$ instead of $\d_{(x_i,\l_i)}$ and $\varphi_{(x_i,\l_i)}$ respectively.

Let  $(\a,\l,x,0)\in M_\e(z)$.  The proof goes along with  the ideas introduced in \cite{BLR}. Once $\ov v$ is defined by Proposition \ref{pv}, then \eqref{AAAZZZ} holds and therefore, by using the Lagrange multiplier theorem, there exist some constants $A$, $B$ and $C$ such that the equation $(E_v)$ (defined in Proposition  \ref{G1}) holds.  Next, we estimate the numbers $A,B,C$ by taking the scalar product of $(E_v)$ with $\d_i$, $\l_i {\partial \d_i}/{\partial \l_i}$ and  $\l_i^{-1}{\partial \d_i}/{\partial x_i}$  respectively. Thus  we derive  a quasi-diagonal system in the variables $A$, $B$ and $C_i$'s. The right  hand side is given by (using Proposition \ref{devIeps} and the fact that $(\a,\l,x,0)\in M_\e(z)$ and ${\partial \Psi_{\e}}/{\partial v}= \n I_\e(u)$ with $u:= \sum_{i=1}^{m}\a_i \d_i  +\ov{v}$)
\begin{eqnarray}\label{F11}
\langle  \frac{\partial \Psi_{\e}}{\partial v}, \d_i\rangle = O (\e \ln ^2\e) ;\quad  \langle \frac{\partial \Psi_{\e}}{\partial v}, \l_i \frac{\partial \d_i}{\partial \l_i}\rangle = O (\e )  ;\quad  \langle \frac{\partial \Psi_{\e}}{\partial v},\frac{1}{\l_i}\frac{\partial \d_i}{\partial x_i}\rangle =  O (\e \ln ^2\e) .\end{eqnarray}
 Hence we deduce that
\be\label{F09}  A_i=O ( \e\ln^2\e) \quad ; \quad B_i=O (\e) \quad ;\quad  C_i=O ( \e\ln^2\e)  \quad \mbox{ for } i=1,\cdots,m.\ee
Furthermore, since $(\a,\l,x,0) \in M_\e(z)$ then Proposition \ref{pv} implies that $\| \ov{v}\| = O(\e)$ and  therefore the system $((E_{\a_i}),(E_{\l_i}),(E_{x_i}))$ (introduced in Proposition \ref{G1}) is  equivalent to

$$ (S_4): \qquad \begin{cases}
&  \langle \n I_\e(u),\d_i \rangle  =  0  \quad \mbox{ for } i = 1,\cdots,m;\\
&  \langle \n I_\e(u),{\l_i}{\partial\d_i}/{\partial \l_i}\rangle = O ( \e^2 \ln^2\e) \quad \mbox{ for } i = 1,\cdots,m; \\
& \langle \n I_\e(u),{\l_i}^{-1}{\partial\d_i}/{\partial x_i}\rangle = O ( \e^2 \ln^2\e)  \quad \mbox{ for } i = 1,\cdots,m.
\end{cases} $$
\noindent
In the following, we  denote by $u:= \sum_{i=1}^{m}\a_i \d_i \,  +  \, \ov{v}$ and we perform  the following change of variables:
\begin{align}
& \b_i:= 1 - \a_i^{4/(n-2)} K(x_i)\, , \quad i=1,\cdots,m;\label{b}\\
& \frac{1}{\l_i} := \frac{ c_4 K(z) }{ c_3  (\partial K / \partial\nu)(z)} \e (1+ \Lambda_i) \, , \quad i=1,\cdots,m;\label{l}\\
& x_i:= \frac{z+\ov{x}_i}{ | z+\ov{x}_i | } \quad \mbox{with } \ov{x}_i:= \g\,  \e^{(n-2)/n} (\ov{\xi}_i + \xi_i) \in T_z(\partial \mathbb{S}^n_+) \, , \quad i=1,\cdots,m  \label{x}
\end{align}
where $(\ov{\xi}_1, \cdots,\ov{\xi}_m)$ is a non-degenerate critical point of $\mathcal{F}_{z,m}$ and $\g$ is a positive constant satisfying
\be \label{gamma}
 \frac{ c_5 \g }{ K(z)} =  \frac{ 2^{n-2} c_2}{\g^{n-1}} \Big( \frac{ c_4 K(z) }{ c_3  (\partial K / \partial\nu)(z)} \Big)^{n-2}.\ee

\noindent
Next, observe that, since  $ \ov{x}_i \in T_z(\partial \mathbb{S}^n_+)$, it follows that
$$ x_i := \frac{z+\ov{x}_i}{ | z+\ov{x}_i | } = (z+\ov{x}_i)( 1 + | \ov{x}_i |^2 )^{-1/2} = z+\ov{x}_i + O(| \ov{x}_i |^2) = z+\ov{x}_i + O( \e^{2(n-2)/n}).$$
Therefore,  under  this change of variables, there hold:
\begin{align}
& d(x_i, x_j) \sim | x_i - x_j | = | \ov{x}_i - \ov{x}_j + O(| \ov{x}_i |^2 + | \ov{x}_j |^2)  | = O(\e^{(n-2)/n}), \mbox{ for all }\,  i,j=1,\cdots,m; \label{aa12}\\
& \e_{ij} \leq \frac{c}{(\l_i\l_j d( x_i , x_j )^2 ) ^{(n-2)/2}} \leq c \, \e^{2(n-2)/n} = o(\e), \mbox{ for all } \, i \neq j .\label{eps12}
\end{align}
Moreover with this change of variables, the first claim of Proposition \ref{devIeps} can be written as follows:
\be\label{1} \langle \n  I_\e(u),\d_i\rangle= \,\a_i \,S_n \b_i+O\big(\e\ln \e\big). \ee
Furthermore it follows from \eqref{x} that
\be \label{kx}
K(x_i) = K_1(x_i) = K_1(z) + O(\e^{2(n-2)/n})   \quad \mbox{ and } \quad \frac{\partial K}{\partial \nu} (x_i) = \frac{\partial K}{\partial \nu} (z) + O(\e^{(n-2)/n}) . \ee
Hence the second claim of Proposition \ref{devIeps} can be written as:
\be\label{2}  \langle \n I_\e(u),\l_i\frac{\partial\d_i}{\partial \l_i}\rangle = -  \a_i^p  c_4 K(z) \e \Lambda_i  + O(\e^{1+ (n-4)/n}).\ee
Next  we deal with  the third claim of Proposition \ref{devIeps}.  First observe that, using the previous change of variables, one derives:
\begin{align*}
& \n_T K(x_i) = \n K_1 (x_i) = \e^{(n-2)/n} \Big( D^2 K_1 (z) (\ov{\xi}_i,.) +  D^2 K_1 (z) ({\xi}_i,.) \Big) + O(\e^{2(n-2)/n}),\\
& \a_i = K(z)^{(2-n)/4} + O(\e \ln^2\e) \qquad ; \qquad  \a_i^p = K(z)^{-(n+2)/4} + O(\e \ln^2\e) ,  \\
& \frac{\partial \e_{ij}}{\partial x_i} = \frac{n-2}{4} \e_{ij}^{n/(n-2)} \l_i\l_j (x_j -x_i) = 2^{n-2}(n-2) \frac{ x_j - x_i } { | x_i - x_j |^n } \frac{1}{ (\l_i \l_j)^{(n-2)/2} }+ O(\e ^{1+2/n})\end{align*}
(by using the fact that : $2(1-\cos d(x_i,x_j)) = | x_i - x_j |^2$). Furthermore there holds:

\begin{align*}
& \frac{ x_j - x_i } { | x_i - x_j |^n } = \frac{ \g \e^{(n-2)/n} [ ( \ov{\xi}_j - \ov{\xi}_i) +( {\xi}_j - {\xi}_i ) + O(\e ^{(n-2)/n}) ] }{\g^n \e^{n-2}\big( | \ov{\xi}_j - \ov{\xi}_i|^2 +2 \langle \ov{\xi}_j - \ov{\xi}_i,  {\xi}_j -{\xi}_i\rangle + | {\xi}_j - {\xi}_i |^2 + O(\e ^{(n-2)/n}) \big)^{n/2}}\\
&\qquad  \qquad\quad   =  \frac{ \e^{(n-2)/n}}{ \g^{n-1} \e^{n-2}} \Big( \frac{ \ov{\xi}_j - \ov{\xi}_i}{ | \ov{\xi}_j - \ov{\xi}_i|^n} - n \frac{ \ov{\xi}_j - \ov{\xi}_i}{ | \ov{\xi}_j - \ov{\xi}_1|^{n+2}} \langle \ov{\xi}_j - \ov{\xi}_i,  {\xi}_j -{\xi}_1\rangle + \frac{ {\xi}_j - {\xi}_i}{ | \ov{\xi}_j - \ov{\xi}_i|^n} + O(| \xi ^2 | + \e ^{\frac{n-2}{n}}) \Big).
\end{align*}
Hence the third claim of Proposition \ref{devIeps} becomes
\begin{align*}\ \langle \n I_\e(u), & \frac{\partial\d_i}{\partial x_i}\rangle_{\lfloor T_{x_i} (\partial \mathbb{S}^n_+)} =   - \frac{\e^{(n-2)/n}}{K(z)^{(n-2)/4}} \Big\{  \frac{c_5 \g }{K(z)} \big( D^2 K_1 (z) (\ov{\xi}_i, . ) + D^2 K_1 (z) ({\xi}_i, . ) \big) \\
& + (n-2)\frac{ 2^{n-2} c_2 }{\g^{n-1}} \Big( \frac{c_4 K(z)}{c_3 (\partial K/\partial \nu)(z)}\Big)^{n-2} \Big(  \sum_{ j\neq i} \frac{ \ov{\xi}_j - \ov{\xi}_i}{ | \ov{\xi}_j - \ov{\xi}_i|^n} - n \frac{ \ov{\xi}_j - \ov{\xi}_i}{ | \ov{\xi}_j - \ov{\xi}_i|^{n+2}} \langle \ov{\xi}_j - \ov{\xi}_i,  {\xi}_j -{\xi}_i\rangle \\
&  + \frac{ {\xi}_j - {\xi}_i}{ | \ov{\xi}_j - \ov{\xi}_i|^n} + \frac{n-2}{2} (\Lambda_i + \Lambda_j) \frac{ \ov{\xi}_j - \ov{\xi}_i}{ | \ov{\xi}_j - \ov{\xi}_i|^n}
+ O(| \xi ^2 | + | \Lambda |^2) \Big)\Big\} + O(\e \ln^2\e).
\end{align*}

Using the fact that $(\ov{\xi}_1, \cdots, \ov{\xi}_m)$ is a critical point of  $ \mathcal{F}_{z,m}$ and the value of $\g$ (see \eqref{gamma}), we obtain
\begin{align}  \langle \n I_\e(u), & \frac{\partial\d_i}{\partial x_i}\rangle_{\lfloor T_{x_i}( \partial \mathbb{S}^n_+)} =   - \frac{\e^{(n-2)/n}}{K(z)^{(n-2)/4}}  \frac{c_5 \g }{K(z)} \Big\{   D^2 K_1 (z) ({\xi}_i, . ) \nonumber \\
& + (n-2)  \Big(   - n  \, \sum_{j \neq i }\frac{ \ov{\xi}_j - \ov{\xi}_i}{ | \ov{\xi}_j - \ov{\xi}_i|^{n+2}} \langle \ov{\xi}_j - \ov{\xi}_i,  {\xi}_j -{\xi}_i\rangle \nonumber\\
&  + \frac{ {\xi}_j - {\xi}_i}{ | \ov{\xi}_j - \ov{\xi}_i|^n} + \frac{n-2}{2} (\Lambda_i + \Lambda_j) \frac{ \ov{\xi}_j - \ov{\xi}_i}{ | \ov{\xi}_j - \ov{\xi}_i|^n}
+ O(| \xi ^2 | + | \Lambda |^2 + \e^{2/n} \ln^2\e) \Big)\Big\} \nonumber\\
& =  - \frac{\e^{(n-2)/n}}{K(z)^{(n-2)/4}}  \frac{c_5 \g }{K(z)} \Big\{    D^2  \mathcal{F}_{z,m} (\ov{\xi}_1, \cdots,  \ov{\xi}_k)( ({\xi}_1, \cdots, \xi_k),.)  \label{3}\\
& + \frac{(n-2)^2}{2}  \sum_{i \neq j}  (\Lambda_i + \Lambda_j) \frac{ \ov{\xi}_j - \ov{\xi}_i}{ | \ov{\xi}_j - \ov{\xi}_i|^n}  + O(| \xi ^2 | + | \Lambda |^2 + \e^{2/n} \ln^2\e) \Big\} .\nonumber
\end{align}
Now using \eqref{1}, \eqref{2} and \eqref{3}, the system $(S_4)$ is equivalent to
\be\label{S'} (S') \, \begin{cases}
 \b_i = O ( \e\ln\e) \, ,  \quad \mbox{ for } i = 1,\cdots, m; \\
 \Lambda_i =   O ( \e^{(n-4)/n} ) \, ,  \quad \mbox{ for } i = 1, \cdots,m \\
D^2  \mathcal{F}_{z,m} (\ov{\xi}_1, \cdots, \ov{\xi}_m) (({\xi}_1, \cdots, \xi_m) ,.) = O (| \xi ^2 | + | \Lambda |^2 + \e^{\frac{2}{n}} \ln^2\e+  \e^{\frac{n-4}{n} }) \, , \mbox{ for } i \leq m.
\end{cases} \ee
Since $(\ov{\xi}_1, \cdots, \ov{\xi}_m)$ is a non-degenerate critical point of $ \mathcal{F}_{z,m}$, we derive that the system $(S')$ has a solution $(\beta_1^\e, \cdots,\beta_m^\e, \Lambda_1^\e, \cdots, \Lambda_m^\e,\xi_1^\e, \cdots, \xi_m^\e)$ and therefore 
$$u_\e := \sum_{i=1}^m \a_i^\e \d_{x_i^\e, \l_i^\e} + \ov{v}_\e $$ (see \eqref{b}-\eqref{x}) is a solution of $(\mathcal{P}_\e)$ which blows up at the same point $z$, as $\e \to 0$.

Finally, we  claim that this constructed solution is positive. In fact, writing $u_\e := u_\e^+ - u_\e^-$, it is easy to get that $| u_\e ^- |_{L^{2n/(n-2)} }\leq | \ov{v}_\e  |_{L^{2n/(n-2)}}=o_\e(1)$. Multiplying Equation $(\mathcal{P}_\e)$ by $u_\e^-$ and integrating over $\mathbb{S}^n_+$ and using the Holder's inequality, we get
$$ \| u_\e ^- \|^2 \leq c | u_\e ^- |^{2n/(n-2)}_{L^{2n/(n-2)} }\leq c \| u_\e ^- \|^{2n/(n-2)}.$$
Hence, either $u_\e^- =0$ or $\| u_\e ^- \| \geq c$ for some constant $c > 0$. The second alternative contradicts the smallness of  $| u_\e ^- |_{L^{2n/(n-2)} }$. Hence we derive that $u_\e^- =0$. Therefore $u_\e > 0$ by using the  maximum principal.\\
Regarding the energy of the solution, observe that, since the interactions between the bubbles are very small, for each $ r > 0$, it holds
\begin{align*}
\int _{B_r (z) } K  u_\e ^{2n/(n-2)} & = \sum_{i=1}^m \int _{B_r (z) } K  (\a_i^\e \d_{x_i^\e, \l_i^\e})  ^{2n/(n-2)}  + o_\e (1) \\
&  = \sum_{i=1}^m (\a_i^\e )  ^{2n/(n-2)}  K( x_i^\e ) \int _{B_r (z) }  ( \d_{x_i^\e, \l_i^\e})  ^{2n/(n-2)}  + o_\e 1) \\
& =  m \frac{1} { K(z)^{(n-2)/2} } S_n  + o_\e (1).
\end{align*}

Concerning the second part of the theorem, which consists of providing a  lower bound  for the Morse index of the constructed  blowing up solutions, we  recall that the Morse index of $I_\e$ at the solution $u_\e$ is defined as the  dimension of the  space  $E$'s such that $$ D^2 I_\e (u_\e) ( \psi , \psi ) < 0 \quad \mbox{ for each } \psi \in E.$$
Next, for  $u_\e:= \sum_{i=1}^{m } \a_{i}^\e \d_{a_i ^\e, \l_i ^\e} + \ov{v}_\e$, observe that
$$ D^2 I_\e (u_\e) ( \d_{a_i ^\e, \l_i ^\e} , \d_{a_i ^\e, \l_i ^\e} ) = \| \d_{a_i ^\e, \l_i ^\e} \|^2 - (p-\e) \int_{\mathbb{S}^n_+} K u_\e ^{p-\e -1} \d_{a_i ^\e, \l_i ^\e}^2.$$
Furthermore notice that, using \eqref{zza1}, we have  that
$$ \int_{\mathbb{S}^n_+} K u_\e ^{p-\e -1} \d_{a_i ^\e, \l_i ^\e}^2 = \a_i ^ {p-1-\e} \int_{\mathbb{S}^n_+} K  \d_{a_i ^\e, \l_i ^\e}^{p+1-\e} + o(1)= \a_i ^ {p-1-\e} K(a_i)\l_i^{-\e(n-2)/2}  S_n +o(1) .$$
Hence, using \eqref{az2}, it follows that $ D^2 I_\e (u_\e) ( \d_{a_i ^\e, \l_i ^\e} , \d_{a_i ^\e, \l_i ^\e} ) \leq S_n (1-p +o(1)) $ for each $i =1, \cdots,m$. In addition, since we have $ \langle
\d_{a_i ^\e, \l_i ^\e} , \d_{a_j ^\e, \l_j ^\e} \rangle =o(1)$  for each $i \neq j$ we derive that the vector space $E:= span\{ \d_{a_1 ^\e, \l_1 ^\e}, \cdots,\d_{a_m ^\e, \l_m ^\e} \}$ is of dimension $m$ and it satisfies that $ D^2 I_\e (u_\e) ( \psi , \psi ) < 0 $ for each $\psi \in E$. Thus we have  that the Morse index is larger   or equal to $m$.

The proof of Theorem \ref{t:11} is thereby complete.

\subsection*{ Proof of Corollary \ref{t:01} }

Using Proposition \ref{pF}, we derive that under the assumption of Theorem \ref{t:01}, the function $\mathcal{F}_{z,2}$ admits a non degenerate critical point. Hence the theorem follows from Theorem \ref{t:11} by taking $m=2$.

\subsection*{ Proof of Theorem \ref{t:13} }

We  follow closely  the proof of Theorem \ref{t:11} but since  we have to change the set $M_\e$, some changes will be needed. In the following proof we will consider  the case where $m \geq 1$ and $\ell \geq 1$, which     corresponds to the most  complicated  situation. Indeed  to recover the other situations, it suffices to remove the corresponding variables.\\
As assumed  in the theorem, let $z_1,\cdots,z_m$ be $m$ critical points of $K_1$ satisfying $\partial K/\partial \nu (z_i)> 0$ and let $y_{m+1}, \cdots, y_{m+\ell}$ be $\ell$ critical points of $K$ satisfying $ \D K(y_j) < 0$. In this proof $M_\e$ will be  defined by:
\begin{align*}  M_\e :=  \{ & \mathcal{M} :=  (\a,\l,x,v)\in(\R_+)^N \times(\R_+)^N\times (\partial \mathbb{S}^n_+)^m \times (\mathbb{S}^n_+)^{\ell } \times H^1(\mathbb{S}^n_+): \\
&   v\mbox{ satisfies \eqref{V0}} \, ; \, \,  |{\a_i^{4/(n-2)}K(x_i)} - 1 | < \e \ln^2\e\, \forall \, i; \, \, C^{-1} \e \leq \l_i^{-1} \leq  C \e \, \forall i \leq m\, ;\\
&  C^{-1} \e \leq \l_j^{-2} \leq  C \e \, \forall j > m \, ; \, \,  d( x_i ,  z_i )  \leq \mu \,\, \forall \, i\leq m \, ; \, \,   d(x_j ,  y_j )  \leq \mu \,\, \forall \, j > m \},
\end{align*}
where $N:= m+\ell$,  $C$ is a large positive constant and  $\mu$ is a small positive one.\\
We notice that, for $(\a,\l,x,0)\in M_\e$, one has that  $ d( x_k , x_r )\geq c >0$ for $k \neq r$ and therefore, using \eqref{epsilon} we derive that $  \e_{kr} = O (\e ^{(n-2)/2})$.
These two properties are the main difference comparing with the proof of Theorem \ref{t:11}.\\
As in the proof of Theorem \ref{t:11}, we define the function $\Psi_\e$ (introduced in \eqref{F10}) and we  note that Proposition \ref{G1} holds with $N=m+\ell$  which implies the existence of  constants $A$, $B$ and $C$ such that the equation $(E_v)$ holds.    In addition, in the same way, as before, \eqref{F11} holds and the constants $A$, $B$ and $C_i$'s satisfy \eqref{F09}.

Therefore the system $((E_{\a_i}),((E_{\l_i}),(E_{x_i}))$ (introduced in Proposition \ref{G1}) is  equivalent to
$$ (S_1): \qquad \begin{cases}
&  \langle \n I_\e(u),\varphi_i \rangle  =  0  \quad \mbox{ for } i = 1,\cdots,N;\\
&  \langle \n I_\e(u),{\l_i}{\partial\varphi_i}/{\partial \l_i}\rangle = O ( \e^2 \ln^2\e) \quad \mbox{ for } i = 1,\cdots,N; \\
& \langle \n I_\e(u),{\l_i}^{-1}{\partial\varphi_i}/{\partial x_i}\rangle = O ( \e^2 \ln^2\e)  \quad \mbox{ for } i = 1,\cdots,N.
\end{cases} $$
\noindent
In the following, we  denote by $u:= \sum_{i=1}^{N}\a_i \varphi_i \,  +  \, \ov{v}$ and we perform  the following change of variables:
\begin{align}
& \b_i:= 1 - \a_i^{4/(n-2)} K(x_i)\, , \quad i=1,\cdots,N;\label{b1}\\
& \frac{1}{\l_i} := \frac{ c_4 K(z_i) }{ c_3  (\partial K / \partial\nu)(z_i)} \e (1+ \Lambda_i) \, , \quad i=1,\cdots,m;\label{l1}\\ 
& \frac{1}{\l_j ^2} := - \frac{ c_4 K(y_j) }{ c_2  \D K (y_j)} \e (1+ \Lambda_j) \, , \quad j = m+1,\cdots,N;\label{l11}\\
& x_i:= \frac{z_i+\xi_i}{| z_i+\xi_i | } \, , \quad i=1,\cdots,m \quad ; \qquad   x_j := \frac{ y_j + \xi_j }{ | y_j+\xi_j | } \, , \quad j = m+1,\cdots, N . \label{x1}
\end{align}
We notice that the main difference between this proof and the one of Theorem \ref{t:11} is the definition of  the concentration points $x_i$'s given by  \eqref{x}  respectively by \eqref{x1}.

Moreover with this change of variables,  for $i \leq m$ (that is $x_i$ is a boundary concentration point),  Proposition \ref{devIeps} implies that \eqref{1} and \eqref{2} hold. Concerning the concentration point, since in this case we have $d(x_k, x_r )  \geq c > 0 $ for each $ k \neq r$, using the third assertion of Proposition \ref{devIeps}, it follows that the counterpart of \eqref{3} becomes
\begin{align}
\langle \n I_\e(u),  \frac{\partial\d_i}{\partial x_i}\rangle_{\lfloor T_{x_i} (\partial \mathbb{S}^n_+)}  & = - c_5 \a _i ^p \l_i^{-\e(n-2)/2} { \n _T K(x_i) }+ O(\e^{1/2}) \nonumber\\
& = - c_5 \a _i ^p  D^2 K_1 (z_i) ( \xi_i,.) + O(\e^{1/2} + | \xi_i |^2) , \qquad \forall \, \, i \leq m.\label{3*3}
\end{align}
Now we will focus on the indices $ j \geq m+1$. This part does not exist in the proof of Theorem \ref{t:11} but the argument is still the same. Using Proposition \ref{devIeps2} and the previous change of variables, it follows that:
\begin{align}
& \langle \n I_\e(u),  \varphi_j \rangle  = 2 \, \a_j S_n\,  \b_j + O\big(\e\ln \e\big)\qquad \forall \, \, m+1 \leq j \leq N, \label{1*1} \\
& \langle \n I_\e(u),  \l_j \frac{\partial\varphi_j }{\partial \l_j }\rangle = - 2\, \e\,  c_4\,  K( y_j ) \Lambda_j  + O(\e^{3/2} +\e  | \xi_j |^2) \qquad \forall  \, \,  m+1 \leq j \leq N, \label{2*2}  \\
& \langle \n I_\e(u),  \frac{\partial\varphi_j }{\partial x_j }\rangle = -2 c_5 \a_j ^p D^2 K ( y_j ) ( \xi_j , .) + O(\e ^{1/2} + | \xi_j |^2) \qquad \forall  \, \, m+1 \leq j \leq N.\label{3**3}
\end{align}
Now using \eqref{1}, \eqref{2}, \eqref{3*3}--\eqref{3**3}, the system $(S_1)$ is equivalent to
\be\label{S'1} (S'_1) \, \begin{cases}
 \b_i = O ( \e\ln\e) \, ,  \quad \mbox{ for } i = 1,\cdots, N; \\
 \Lambda_i =   O ( \e^{1/5} ) \, ,  \quad \mbox{ for } i = 1, \cdots,N \\
D^2 K_1 (z_i) ( \xi_i ,.) = O (| \xi ^2 |  + \e^{{2}/{n}} \ln^2\e+  \e^{{1}/{5} }) \, , \mbox{ for } i \leq m,\\
D^2 K ( y_j ) ( \xi_j ,.) = O (| \xi ^2 |  + \e^{{2}/{n}} \ln^2\e+  \e^{{1}/{5} }) \, , \mbox{ for } j \geq m+1 .
\end{cases} \ee
We derive that the system $(S'_1)$ has a solution $(\beta_1^\e, \cdots,\beta_N^\e, \Lambda_1^\e, \cdots, \Lambda_N^\e,\xi_1^\e, \cdots, \xi_N^\e)$ and therefore $u_\e := \sum \a_i^\e \varphi_{x_i^\e, \l_i^\e} + \ov{v}_\e $ (see \eqref{b1}-\eqref{x1}) is a solution of $(\mathcal{P}_\e)$ which blows up at $N=m+\ell$ simple blow up points $\{z_1,\cdots,z_m,y_{m+1},\cdots,y_{m+\ell}\}$, as $\e \to 0$.

Finally, arguing as in the end of the proof of Theorem \ref{t:11}, we derive  that the  constructed function $u_\e $ is positive.
The proof of Theorem \ref{t:13} is thereby complete.

\subsection*{ Proof of Theorems \ref{t:18}  }

The proof follows the proof  of Theorem \ref{t:13} step by step with some changes. In fact, the set $M_\e$ will be
$$ M'_\e:= \{ (\a_0,\a,\l,x,v) \in \R\times M_\e : v \mbox{ satisfies  \eqref{wV0} } \, ; \, | \a_0^{4/(n-2)}-1 | \leq \e \ln^2(\e) \}.$$
The function $\Psi_\e$ will be changed as
\begin{eqnarray}\label{0*01}  \Psi_{\e}:M'_{\e}\,\rightarrow \R ;\quad \mathcal{M}=(\a_0,\a,\l,x,v)\mapsto I_\e\big( \a_0 \o + \sum_{i=1}^{N} \a_i \varphi_{(x_i,\l_i)} \, + \, v\big). \end{eqnarray}
In Proposition \ref{G1}, we will have another equation $(E_{\a_0})$ corresponding to the variable $\a_0$ which gives another equation in the system $(S_1)$:
\be\label{0*02} \langle \n I_\e(u) , \o \rangle = 0.\ee
 Furthermore, since $v$ satisfies \eqref{wV0} hence it is also orthogonal to the function $\o$. Then, once $\ov{v}$ is found by Proposition \ref{wpv}, the equation \eqref{AAAZZZw} implies the existence  of some constants $A_0$, $A$, $B$ and $C$ such that
\be\label{AAAZZZ1} (E_v ^\o): \qquad \frac{\partial \Psi_{\e}}{\partial v}(\a,\l,x,v)= A_0 \o + \sum_{i=1}^N \biggl(A_i \varphi_i+B_i\l_i \frac{\partial \varphi_i}{\partial \l_i}+\sum_{j=1}^{n} C_{ij}\frac{1}{\l_i}\frac{\partial \varphi_i}{\partial x_i^j}\biggr) \ee
which is the new equation $(E_v)$ when $\o \neq 0$.

In \eqref{b1}, we will add another equation which is $\b_0 := 1 - \a_0^{4/(n-2)}$.\\
The sequel of the proof is unchanged but we will use the corresponding formulae from Subsection 3.1.

\subsection*{ Proof of Theorem \ref{t:12}}

For $\ell =0$ (that is  when there is  only boundary blow up points), the proof follows the proof  of Theorem \ref{t:11} step by step getting  more equations in the system $(E_{\a_i}, E_{\l_i}, E_{a_i}, E_v)$. Indeed  in this proof we modify the set $M_\e(z)$ as follows:
\begin{align} M_\e(z_1,\cdots,z_m):= \{ ({ \a}^1, \cdots, & { \a}^m, { \l}^1,\cdots,{ \l}^m,  { a}^1,\cdots,{ a}^ m,v): \label{0*00}\\
& (\a^i, \l^i , a^ i , 0) \in M_\e(z_i) \, \, \forall \, i  \mbox{ and } v \mbox{ satisfies \eqref{V0}}\}. \nonumber \end{align}
Note that, in this proof we have $\a ^i = ( \a ^i _1, \cdots , \a^i_{q_i})$, $\l^i = ( \l^i _1, \cdots , \l^i_{q_i})$ and $a^i = ( a^i _1, \cdots , a^i_{q_i})$ for each $i=1, \cdots, m$ and $ d( a^i _k , z_i ) \to 0$ for each $k$. Furthermore,  for $a^i_j$ close to $z_i$ and $a^k_\ell$ close to $z_k$ with $k \neq i$, it follows that $ d( a^i_j , a^k_\ell ) \geq c $  and therefore  the interaction between the two corresponding bubbles is of the order of
$$ O\Big(\frac{1}{( \l_j^i \l_\ell^k)^{(n-2)/2} d(a^i _j , a^k_j ) ^{n-2}}\Big) = O(\e^{n-2})$$ since the rates $\l$'s are of order $1/\e$ (see the definition of $M_\e(z_i)$ defined in \eqref{Meps}). Hence this interaction does not have any contribution in the phenomenon.

Following   the proof  of Theorem \ref{t:11}, the first step  which consists of  finding   and estimating the  infinite dimensional part $\ov{v}$  is already  done in  Proposition \ref{pv}. The second step consists to solve the following system:
 $$ (S'_4 ): \qquad \begin{cases}
&  \langle \n I_\e(u),\d_{a^ i_j, \l ^ i  _j} \rangle  =  0  \quad \mbox{ for } i = 1,\cdots,m \, \mbox { and } \forall \, \, j = 1 , \cdots , q_i ; \\
&  \langle \n I_\e(u),{\l^i _j}{\partial\d_{a^ i_j, \l ^ i  _j}}/{\partial \l^i _j}\rangle = O ( \e^2 \ln^2\e) \quad \mbox{ for } i = 1,\cdots,m \, \mbox { and } \forall \, \, j = 1 , \cdots , q_i ; \\
& \langle \n I_\e(u),({\l^i _j})^{-1}{\partial\d_{a^ i_j, \l ^ i  _j}}/{\partial a^i _j}\rangle = O ( \e^2 \ln^2\e)  \quad \mbox{ for } i = 1,\cdots,m \, \mbox { and } \forall \, \, j = 1 , \cdots , q_i .
\end{cases} $$
We remark that this system is similar to $(S_4)$ introduced in the proof of Theorem \ref{t:11} but it contains more equations. To study this system, we introduce the following change of variables :
\begin{align*}
& \b^i_j := 1 - (\a^i_j) ^{4/(n-2)} K(a^i _j )\, , \quad i=1,\cdots,m  \, ; \, \, \forall \, \, j = 1 , \cdots , q_i ;\\
& \frac{1}{\l^i_j } := \frac{ c_4 K(z_i) }{ c_3  (\partial K / \partial\nu)(z_i)} \e (1+ \Lambda ^i_j ) \, , \quad i=1,\cdots,m  \, ; \, \,  \forall \, \, j = 1 , \cdots , q_i ; \\
& a^i_j := \frac{z+\ov{x}^i_j }{ | z+\ov{x} ^i_j } \quad \mbox{with } \ov{x}^i_j := \g_i \,  \e^{(n-2)/n} (\ov{\xi}^i_j + \xi ^i_j ) \in T_{z_i}(\partial \mathbb{S}^n_+) \, , \quad i=1,\cdots,m   \,  ; \, \, \forall \, \, j = 1 , \cdots , q_i
\end{align*}
where $(\ov{\xi}^i_1 , \cdots,\ov{\xi} ^i_{q_i} )$ is a non-degenerate critical point of $\mathcal{F}_{z_i,q_i}$ and $\g_i$ is a positive constant satisfying
$$
 \frac{ c_5 \g _i }{ K(z_i)} =  \frac{ 2^{n-2} c_2}{\g_i ^{n-1}} \Big( \frac{ c_4 K(z_i) }{ c_3  (\partial K / \partial\nu)(z_i)} \Big)^{n-2}.$$
Following the computations performed  in the proof of Theorem \ref{t:11} we derive that the system $(S'_4)$ is equivalent to the following one:
$$ (S'') \, \begin{cases}
 \b^ i_j = O ( \e\ln\e) \, ,  \quad \mbox{ for } i = 1,\cdots, m  \, ; \, \, \forall \, \, j = 1 , \cdots , q_i  ; \\
 \Lambda^i_j  =   O ( \e^{(n-4)/n} ) \, ,  \quad \mbox{ for } i = 1, \cdots,m  \, ; \, \, \forall \, \, j = 1 , \cdots , q_i ;\\
D^2  \mathcal{F}_{z_i,q_i} (\ov{\xi}^i_1, \cdots, \ov{\xi}^i_{q_i}) (({\xi}^i_1, \cdots, \xi^i_{q_i}) ,.) = O (| \xi ^2 | + | \Lambda |^2 + \e^{\frac{2}{n}} \ln^2\e+  \e^{\frac{n-4}{n} }) \, , \mbox{ for } i \leq m ,
\end{cases} $$
where $ \xi := (\xi ^1, \cdots, \xi ^m)$,  $  \Lambda := (\Lambda^1, \cdots, \Lambda ^m)$, $\xi ^ i := ( \xi ^i_1, \cdots, \xi ^i _{q_i})$ and $\Lambda ^ i := ( \Lambda ^i_1, \cdots, \Lambda ^i _{q_i})$.

We recall that, for each $i= 1, \cdots,m$, we have $(\ov{\xi}^i _{1}, \cdots, \ov{\xi}^i _{q_i})$ is a non-degenerate critical point of $ \mathcal{F}_{z_i ,q_i}$. Thus, by using  Browder-fixed point theorem, we derive that the previous system $(S'')$ has a solution $(\beta^{1, \e} , \cdots,\beta^{m, \e}, \Lambda ^{1, \e}, \cdots, \Lambda^{m, \e }, \xi^{1, \e}, \cdots, \xi^{ m, \e} )$ and therefore $u_\e := \sum_{i,j} \a_j^{i, \e} \d_{a_j^{i, \e}, \l _j^{i, \e}} + \ov{v}_\e $  is a solution of $(\mathcal{P}_\e)$ which blows up at the points $ ( z_1, \cdots,z_m)$, as $\e \to 0$.

 As in the end of the proof of Theorem \ref{t:11}, this function is positive  and the proof of Theorem \ref{t:12} is thereby completed in the case where $\ell = 0$.\\
Concerning the case where $\ell \geq 1$, we can combine the argument used  in the proof of Theorem \ref{t:13} for the variables $\a_i$, $a_i$ and $\l_i$ corresponding to the blow up point $y_i$ and the previous argument to derive the result (as in the proof of Theorem \ref{t:13}). In fact, the new system looks as  $(S'', S''_{int})$ where $S''_{int}$ will contain the indices $i > \sum _{j \leq m} q_j$ as in the proof of Theorem \ref{t:13}.\\
Hence the proof of Theorem \ref{t:12} is thereby completed.

 \section*{ Proof of   Theorem \ref{t:15} }

 Let $\o$ be a non degenerate solution of $(\mathcal{P})$. In a first step we prove the  first statement in Theorem \ref{t:15}\\
 {\it Proof of statement  $a)$}:
 Let $z_1,\cdots,z_m$ be $m$ critical points of $K_1$ with $\partial K/\partial (z_i) > 0$.\\
  Once again we follow the proof   of Theorem \ref{t:12}. Actually  compared  with the proof of Theorem \ref{t:12}, we will have to add  another variable $\a_0$ and the function $u$ looks like
 $$u= \a_0 \o + \sum \a_i \d_i + v.$$ The new set $M_\e$ will be defined by
 $$ M_\e:= \{ (\a_0, \a,\l,x,v): |\a_0^{4/(n-2)} - 1 | \leq \e \ln^2 \e \,  ; \, (\a,\l,x,v)\in M_\e(z_1,\cdots ,z_m)\, ; \, v \mbox{ satisfies \eqref{wV0}}\}$$
 where $M_\e(z_1,\cdots,z_m)$ is introduced in \eqref{0*00}.
 Furthermore the function $\Psi_\e$ will be as introduced in \eqref{0*01}. Hence with respect to Proposition \ref{G1} we will have another equation $(E_{\a_0})$ corresponding to the variable $\a_0$ which is defined by
 $$  (E_{\a_0}) \qquad \frac{\partial \Psi_{\e}}{\partial \a_0}(\a_0, \a,\l,x,v)=0.$$
 Furthermore, once $\ov{v}$ is defined by Proposition \ref{wpv}, by the Lagrange multiplier theorem, we derive from   \eqref{AAAZZZw}    the existence of some constants $A_0$, $A$, $B$ and $C$ such that
$$ (E_v ^\o): \qquad \frac{\partial \Psi_{\e}}{\partial v}(\a,\l,x,v)= A_0 \o + \sum_{i=1}^N \biggl(A_i \varphi_i+B_i\l_i \frac{\partial \varphi_i}{\partial \l_i}+\sum_{j=1}^{n} C_{ij}\frac{1}{\l_i}\frac{\partial \varphi_i}{\partial x_i^j}\biggr). $$
Once $E_v^\o$ is satisfied, it remains to solve the following system

 $$ (S''_4 ): \qquad \begin{cases}
 &  \langle \n I_\e(u),\o \rangle  =  0 ; \\
&  \langle \n I_\e(u),\d_{a^ i_j, \l ^ i  _j} \rangle  =  0  \quad \mbox{ for } i = 1,\cdots,m \, \mbox { and } \forall \, \, j = 1 , \cdots , q_i ; \\
&  \langle \n I_\e(u),{\l^i _j}{\partial\d_{a^ i_j, \l ^ i  _j}}/{\partial \l^i _j}\rangle = O ( \e^2 \ln^2\e) \quad \mbox{ for } i = 1,\cdots,m \, \mbox { and } \forall \, \, j = 1 , \cdots , q_i ; \\
& \langle \n I_\e(u),({\l^i _j})^{-1}{\partial\d_{a^ i_j, \l ^ i  _j}}/{\partial a^i _j}\rangle = O ( \e^2 \ln^2\e)  \quad \mbox{ for } i = 1,\cdots,m \, \mbox { and } \forall \, \, j = 1 , \cdots , q_i .
\end{cases} $$
We notice that the first equation is dealt with  in Proposition  \ref{devIeps4}, while the other ones are considered  in Proposition  \ref{devIeps3}. Moreover observe  that these propositions are similar to Propositions  \ref{devIeps} and  \ref{devIeps2} and the main difference comes from  the contribution of the function $\o$. Furthermore  notice that, in the first and the third Claims of Proposition  \ref{devIeps3} and in the second and the fourth Claims of Proposition  \ref{devIeps4} , this contribution is seen as a remainder term. However, in the third Claim of Proposition  \ref{devIeps4}, this contribution is of the order of
$$  \frac{ \o(a_i) }{\l_i ^{(n-2)/2}} $$   and might be  a principal term, depending  on the dimension $n$.

 Recall that in Claim $a)$ of Theorem \ref{t:15} we have only boundary blow up points. In this case, in the second claim of Proposition  \ref{devIeps3}, we have $ \e$ and  $ \frac{1}{\l } $ as  principal terms and therefore the contribution of $\o$ (which is in this case of the order of $O( \e^{(n-2)/2})$)  is negligible  with respect to the principal parts. Hence this  contribution  comes  as a remainder term.  Hence as in the proof of Theorem \ref{t:11}, the system $(S''_4)$ is equivalent to the following one:
 $$ (S''') \, \begin{cases}
 \b_0 = O ( \e\ln\e) ; \\
 \b^ i_j = O ( \e\ln\e) \, ,  \quad \mbox{ for } i = 1,\cdots, m  \, ; \, \, \forall \, \, j = 1 , \cdots , q_i  ; \\
 \Lambda^i_j  =   O ( \e^{(n-4)/n} ) \, ,  \quad \mbox{ for } i = 1, \cdots,m  \, ; \, \, \forall \, \, j = 1 , \cdots , q_i ;\\
D^2  \mathcal{F}_{z_i,q_i} (\ov{\xi}^i_1, \cdots, \ov{\xi}^i_{q_i}) (({\xi}^i_1, \cdots, \xi^i_{q_i}) ,.) = O (| \xi ^2 | + | \Lambda |^2 + \e^{\frac{2}{n}} \ln^2\e+  \e^{\frac{n-4}{n} }) \, , \mbox{ for } i \leq m .
\end{cases} $$
We recall that by assumption, for each $i= 1, \cdots,m$, we have  that $(\ov{\xi}^i _{1}, \cdots, \ov{\xi}^i _{q_i})$ is a non-degenerate critical point of $ \mathcal{F}_{z_i ,q_i}$. Thus, by using  Browder-fixed point theorem, we derive the existence of a solution of the previous system $(S''')$ and therefore the existence of a solution $u_\e$ follows.

\medspace

 {\it Proof of statement  $b)$}: In this part, we have at least one interior blow up point $y_i$. We remark that, our constructions are based essentially on the expansions of the gradient of $I_\e$. In this case, we need to use Subsection 3.1 (since $\o \neq 0$)  and precisely, we will use Proposition \ref{devIeps4} (since we have an interior blow up point). In the third assertion of this proposition, for $a_i$ close to $y_i$, there are three principal terms which are : $  c\, \e >0$ , $   c \frac{\D K(a_i)}{\l_i^2} < 0$ and $c \, \frac{\o(a_i)}{\l_i^{(n-2)/2}} > 0$. Note that,

\begin{itemize}
\item for $n = 5$, we derive that $(n-2)/2 = 3/2$ and therefore in this case the second term (which is $   c \frac{\D K(a_i)}{\l_i^2} < 0$) will be small with respect to the third one (which is $c \, \frac{\o(a_i)}{\l_i^{3/2}} > 0$). Hence we obtain only two principal terms having the same sign. Precisely, in the system to solve there exists the following equation
$$ (E_{\l_i}) \qquad \frac{2 c_4 \a_i ^p } { \l_i ^{ \e(n-2)/2} } K(a_i) \e (1 + o(1)  ) + 2 \ov{c}_6 \a_0 \frac{\o (a_i) }{\l_i ^{3/2}} (1+ o(1) ) = 0.$$

 and therefore the corresponding system is not solvable.
\item for $n=6$, the second and the third terms are  of the same order. In this case our argument fails since  we cannot compare $\o(y_i)$ with $\D K(y_i)$.
\item Concerning the case $n \geq 7$, in this case we have $(n-2)/2 > 2$ and therefore the third term will be small with respect to the second one  (that is $ \frac{\o(a_i)}{\l_i^{(n-2)/2}} = o ( \frac{\D K(a_i)}{\l_i^2} ) $)
and it will be seen as a remainder term. Thus the estimates of the claims 2-4 of Proposition \ref{devIeps4} become as the ones in Proposition \ref{devIeps2} and therefore
our argument holds by adding another equation to the proof of Theorem \ref{t:12} (the case $\ell \geq 1$).
 \end{itemize}
 The proof is thereby complete

\section{ Appendix}

In this section we provide various pointwise and integral estimates for the bubble. These estimates are used in the expansion of the gradient near the neighborhood at infinity.

 \begin{lem}\label{lowerL2}
 Let $a\in \ov{ \mathbb{S}^n_+}$ and $\l > 0$ be large.\\
  $(i)$ Assume that $\e \ln \l$ is small enough, then it holds
 \begin{align} \label{lower2}\d_{a,\l}^{-\e}(x) = & c_0^{-\e} \l^{-\e (n-2)/2}\Big(1 + \frac{n-2}{2}\, \e \ln (2+(\l^2-1)(1 -\cos d(a,x))) \Big) \\
 & + O\Big(\e ^2 \ln (2+(\l^2-1)(1 -\cos d(a,x)))\Big) \quad \mbox{ for each } y \in \ov{\mathbb{S}^n_+}. \nonumber \end{align}
 $(ii)$ For each $\g > 0$ and each $\b\in [0, n/(n-2))$, it holds
 $$ 0 < \int_{\mathbb{S}^n_+} \d_{a,\l}^{p+1 - \b } (x) \ln^\g \Big(2+(\l^2-1)(1 -\cos d(a,x))\Big)dx = O\Big( \frac{1}{\l^{\b (n-2)/2} }\Big).$$
 $(iii)$ For $\g = 1$, we can be more precise and we get
 \begin{align*}
 & \int_{\mathbb{S}^n_+} \d_{a,\l}^{p+1} \ln \Big(2+(\l^2-1)(1 -\cos d(a,x))\Big)dx  = \begin{cases}
 \ov{c}_{10} + O\Big( \frac{1}{\l^2} \Big) \quad \mbox{ if } a\in \partial \mathbb{S}^n_+, \\
 2 \ov{c}_{10} + O\Big( \frac{1}{\l^2} + \frac{\ln(\l d)}{(\l d )^n} \Big) \quad \mbox{ if } a\in  \mathbb{S}^n_+, \end{cases}\\
 &  \int_{\mathbb{S}^n_+} \d_{a,\l}^{p} \l \frac{\partial \d_{a,\l} }{\partial \l } \ln \Big(2+(\l^2-1)(1 -\cos d(a,x))\Big)dx  = \begin{cases}
- \ov{c}_{11} + O\Big( \frac{1}{\l^2} \Big) \quad \mbox{ if } a\in \partial \mathbb{S}^n_+, \\
 - 2 \ov{c}_{11} + O\Big( \frac{1}{\l^2} + \frac{\ln(\l d)}{(\l d )^n} \Big) \quad \mbox{ if } a\in  \mathbb{S}^n_+, \end{cases}
\end{align*}
 where
 $$  \ov{c}_{10} :=  \int_{\mathbb{R}^n_+} \frac{ c_0^{p+1}}{(1+| z |^2)^n} \ln \Big( 2+ 2 | z |^2 \Big)  \quad \mbox{ and } \quad
  \ov{c}_{11} :=  c_0^{p+1} \frac{n-2}{2} \int_{\mathbb{R}^n_+} \frac{ | z |^2 - 1 }{(1+| z |^2)^{n+1}} \ln \Big( 1+  | z |^2 \Big) > 0.$$
  \end{lem}
  \begin{pf}
  To prove Claim $(i)$, from the definition of $\d_{a,\l}$ (see \eqref{dal}, it follows that

  $$ \d_{a,\l}^{-\e}(x) =  c_0^{-\e} \l^{-\e (n-2)/2} exp\Big(  \frac{n-2}{2}\, \e \ln (2+(\l^2-1)(1 -\cos d(a,x))) \Big).$$
  Observe that $2 \leq 2+(\l^2-1)(1 -\cos d(a,x))) \leq 2\l^2$ and therefore $\e \ln (2+(\l^2-1)(1 -\cos d(a,x))) $ is small uniformly in $x \in  \mathbb{S}^n_+$. Thus Claim $(i)$ follows from the fact that $exp(1+t) = 1+ t + O(t^2)$ for $t$ small.

 First we will focus on the case where $a \in \partial \mathbb{S}^n_+$. For the second claim, let $\g > 0$ and $\b \in [0,1)$, the left inequality holds since the function is positive. Concerning the right one,  using the change of variables $ y = \pi_{-a}(x)$ (without loss of the generality, we can assume that $a=(-1, 0,\cdots,0,)$) where
  \be\label{chvar} \pi_{-a} ^{-1} : \R^n_+ \to  \mathbb{S}^n_+ \, \, ; \quad y := (y_1,\cdots,y_n) \mapsto x := \frac{2}{1+| y |^2} ( \frac{| y |^2 - 1}{2} , y_1, \cdots, y_n).\ee
  Note that, using this change of variables, it holds that
  $$dx = \Big(\frac{2 }{1+ | y |^2} \Big)^n dy \, \, ; \, \, 1 -\cos d(a,x) = \frac{2 | y | ^2}{1+ | y |^2} \mbox{ and } 2+(\l^2-1)(1 -\cos d(a,x) = 2 \frac{1+ \l^2 | y |^2}{1+ | y |^2}.$$
  Hence we get
  \begin{align*} \int_{\mathbb{S}^n_+}  \cdots & = c_0^{p+1 - \b } \int_{\mathbb{R}^n_+} \frac{ \l^{n- \b (n-2)/2}}{(1+\l^2| y |^2)^{n - \b (n-2)/2} } \ln^\g \Big( \frac{2 }{1+ | y |^2}(1+ \l^2 | y |^2) \Big) \Big(\frac{2 }{1+ | y |^2} \Big)^{\b (n-2)/2} dy \\
  & \leq  c \int_{\mathbb{R}^n_+} \frac{ \l^{n- \b (n-2)/2}}{(1+\l^2| y |^2)^{n- \b (n-2)/2}} \ln^\g \Big( 2+ 2\l^2 | y |^2 \Big)dy\\
&   \leq  \frac{c}{\l^{  \b (n-2)/2}} \int_{\mathbb{R}^n_+} \frac{ 1}{(1+| z |^2)^{n- \b (n-2)/2}} \ln^\g \Big( 2+ 2 | z |^2 \Big)dz
  \end{align*}
  (in the last equality, we  used the change of variables $z = \l y$). Thus the second claim follows.\\
  Concerning the last one, taking $\b = 0$, it holds that
\begin{align*} \int_{\mathbb{S}^n_+}  \cdots & = c_0^{p+1} \int_{\mathbb{R}^n_+} \frac{ \l^n}{(1+\l^2| y |^2)^n} \ln \Big( \frac{2 }{1+ | y |^2}(1+ \l^2 | y |^2) \Big)dy \\
    & = c_0^{p+1} \int_{\mathbb{R}^n_+} \frac{ 1}{(1+| z |^2)^n} \ln \Big( 2+ 2 | z |^2 \Big)dz - c_0^{p+1} \int_{\mathbb{R}^n_+} \frac{ \l^n}{(1+ \l^2 | y |^2)^n} \ln \Big( 1+  | y |^2 \Big)dy.
   \end{align*}
  It remains to estimate the last integral. In fact, since : $ \ln(1+t) \leq t$ for each $t \geq 0$, it follows that
  \be\label{app:1}  \int_{\mathbb{R}^n_+} \d^{p+1}_{0,\l} \ln \Big( 1+  | y |^2 \Big)dy  \leq \int_{\mathbb{R}^n_+} \frac{ c_0^{p+1} \l^n }{(1+ \l^2 | y |^2)^n}  | y|^2 dy \leq \frac{ c}{\l^2}.\ee
  The proof of the second claim of $(iii)$ follows exactly by the same way than the previous one.\\
  Now taking $a  \in  \mathbb{S}^n_+$, the same computations done to prove Claim $(ii)$ hold by taking $ \pi_{-a} ^{-1} : \R^n \to  \mathbb{S}^n $. Hence Claim $(ii)$ is proved in the general case that is $a\in \ov{ \mathbb{S}^n_+}$. Concerning Claim $(iii)$, observe that
  $$  \int_{\mathbb{S}^n_+} \d_{a,\l}^{p+1} \ln \Big(2+(\l^2-1)(1 -\cos d(a,x))\Big)dx  =  \int_{\mathbb{S}^n } \cdots + O\Big( \frac{\ln (\l d) }{ (\l d )^n } \Big).$$
  Hence the proof follows following the case $a \in \partial \mathbb{S}^n_+$ by taking $ \pi_{-a} ^{-1} : \R^n \to  \mathbb{S}^n $.
  This ends the proof.
   \end{pf}

\begin{lem}\cite{AB20a}\label{lem:varphi}
For $a\in \partial \mathbb{S}_+^n$, we have $\partial \delta_{a,\l} /\partial \nu = 0$ and therefore $\varphi_{a,\l} = \d_{a,\l}$. For $a\notin \partial \mathbb{S}_+^n$, we have
$$ \d_{a,\l} \leq  \varphi_{a,\l} \leq 2 \d_{a,\l}\, \, ; \quad | \l \partial \varphi_{a,\l} /\partial \l | \leq c \d_{a,\l} \, \, ; \quad  | (1/\l) \partial \varphi_{a,\l} /\partial a^k | \leq c \d_{a,\l}, \leqno{(i)}$$
where $a^k$ denotes the $k$-th component of $a$.
$$   \varphi_{a,\l}  = \d_{a,\l} + c_0 \frac{H(a,.)}{\l^{(n-2)/2}} + f_{a,\l}\qquad \mbox{ where } \leqno{(ii)} $$
$$ | f_{a,\l} | _{\infty} \leq  \frac{c}{(\l d_a)^2} \frac{H(a,.)}{\l^\frac{n-2}{2}} \leq  \frac{c}{\l^\frac{n+2}{2} d_a^{n}} \, \, ; \, \,  | \l \frac{\partial f_{a,\l} }{\partial \l}  | _{\infty} \leq \frac{c}{\l^\frac{n+2}{2} d_a^{n}} \mbox{ and }  | \frac{1}{\l} \frac{\partial f_{a,\l} }{\partial a^k}  | _{\infty} \leq \frac{c}{\l^\frac{n+4}{2} d_a^{n+1}},$$ where $d_a:= d(a,\partial \mathbb{S}^n_+$). Furthermore, it holds that ${H(a,.)}/{\l^{(n-2)/2}} \leq c \d_{a,\l}$.
 \end{lem}

In the following lemma we collect some estimates which are used essentially in Section 3. The first, the second and the third  ones are quoted from  \cite{B1} (see $(E1)$, $(E2)$  page 4 and $(F16)$ page 23) and the proof of the other ones  follows immediately by using some basic computations.
\begin{lem}\label{lem:estimates} Assume that the variables $\l_k$'s are large and the $\e_{ij}$'s are small. Then there hold
\begin{align}
& \int_{\R^n} \d_i^{(n+2)/(n-2)} \d_j = c_2 \e_{ij} + O( \e_{ij}^{n/(n-2)} ) \quad \mbox{ with } c_2:=  \int_{\R^n} \frac{ c_0^{2n/(n-2)}}{(1+| x |^2)^{(n+2)/2}}dx ,  \label{91'} \\
& \int_{\R^n} \d_j^{(n+2)/(n-2)} \l_i \frac{\partial \d_i}{\partial \l_i}  = c_2  \l_i \frac{\partial  \e_{ij}}{\partial \l_i}  + O( \e_{ij}^{n/(n-2)}  \ln \e_{ij}^{-1}) , \label{91''}\\
&  \int_{\R^n} (\d_k \d_j)^{n/(n-2)} \leq c \,  \e_{kj} ^{n/(n-2)}\ln(\e_{kj}^{-1})  \quad \forall k \neq j \, , \label{za1} \\
&  \int_{\R^n \setminus B(a,d) } \d_{a,\l} ^{ \g } \leq c \, \frac{1}{ \l^{n - \g (n-2)/2} (\l d )^{(n-2)\g - n }} , \quad \forall \, \, \frac{n}{n-2} < \g \leq \frac{2n}{n-2} , \label{2***}\\
& c_0 \int_{\R^n} \d_{a,\l} ^p = \frac{c_2}{\l^{(n-2)/2}} , \label{3***} \\
&  p c_0  \int_{\R^n} \d_{a,\l} ^{p-1} \l \frac{\partial \d_{a,\l}}{\partial \l}  = - \frac{n-2}{2}\frac{c_2}{\l^{(n-2)/2}} , \label{4***}\\
& \Big( \int_{ B(a_i, d_i) } \d_i ^{8n/(n^2-4)} \Big)^{(n+2)/(2n)} \leq c \begin{cases}
1/ \l_i ^{(n-2) /2} \mbox{ if } n \leq 5, \\
\ln^{2/3}(\l_i d_i) / \l_i^2 \mbox{ if } n = 6 ,\\
d_i ^{(n-6)/2} / \l_i ^2 \mbox{ if } n \geq 7. \end{cases} \label{1***}
\end{align}

\end{lem}

{\small

\bigskip

 {\small Mohameden Ahmedou,\\
  Mathematisches Institut  der Justus-Liebig-Universit\"at Giessen,  Arndtsrasse 2, D-35392 Giessen, Germany, \\
  Mohameden.Ahmedou@math.uni-giessen.de  }
\bigskip

 {\small Mohamed Ben Ayed,\\
  Department of Mathematics, College of Science, Qassim University, Buraidah 51452, Saudi Arabia, \\
  M.BenAyed@qu.edu.sa \\ \& \\
 Universit\'e de Sfax,  Facult\'e des Sciences de Sfax,  D\'epartement de Math\'ematiques,  Route de Soukra, Sfax, BP. 1171, 3000, Tunisia,  \\ Mohamed.Benayed@fss.rnu.tn}

\end{document}